\documentclass[12 pt,A4]{article}
\usepackage{amsmath}
\usepackage{amsthm}
\usepackage{amssymb}

\newtheorem{prop}{Proposition}
\newtheorem{theorem}{Theorem}
\newcommand{\mm}{Minkowski}
\newtheorem{lemma}{Lemma}

\usepackage{graphicx}
\usepackage{hyperref}
\numberwithin{equation}{section}
\begin{document}
 \baselineskip16pt

\title{\bf On conjectures of Minkowski and Woods for $n=10$}
\author{Leetika Kathuria\\{\small{\em Mehr Chand Mahajan DAV College for Women, Chandigarh, India.}}\\ Madhu Raka{\footnote{\em The research  is supported by CSIR, sanction no. ES/ 21(1042)/17/ EMR-II.}}\\
 {\small{\em Centre for Advanced Study in Mathematics, Panjab University, Chandigarh, India.}}
\date{}}
\maketitle
\begin{abstract}  Let $\mathbb{L}$ be a lattice in $n$-dimensional Euclidean space $\mathbb{R}^n$  reduced in the sense of Korkine and
 Zolotareff and having a basis of the form  $~(A_1,0,0,\cdots$ $,0),$ ~$(a_{2,1},A_2,0,\cdots,0),\cdots,$ $(a_{n,1},a_{n,2},\cdots,a_{n,n-1},A_n)$. A famous conjecture of Woods
 in  Geometry of Numbers asserts that if $A_1A_2\cdots
A_n = 1$ and $A_i\leq A_1$ for each $i$ then any closed sphere in
$\mathbb{R}^n$ of radius $\sqrt{n/4}$ contains a point of $\mathbb{L}.$  Together with a result of C. T. McMullen (2005), the truth of Woods' Conjecture for a fixed $n$, implies the
long standing classical conjecture of Minkowski on product of $n$ non-homogeneous linear forms for that value of $n$. In an earlier paper `Proc. Indian Acad. Sci. (Math. Sci.) Vol. 126, 2016,  501-548' we proved Woods'  Conjecture for $n=9$. In  this paper,  we prove  Woods' Conjecture and  hence Minkowski's Conjecture for $n=10$.
 \end{abstract}

\section{Introduction}

$~~~$Let $\mathbb{L}$ be a lattice in the Euclidean space $\mathbb{R}^n.$ By
 the reduction theory of quadratic forms introduced by Korkine and
 Zolotareff \cite{KZ}, a cartesian co-ordinate system may be chosen in
$\mathbb{R}^n$
 in such a way that $\mathbb{L}$ has a basis of the form
$$~(A_1,0,0,\cdots,0),~(a_{2,1},A_2,0,\cdots,0),\cdots,(a_{n,1},a_{n,2},\cdots,a_{n,n-1},A_n),$$
 where $A_1,A_2,$ $\cdots,A_n$ are all positive and further for each $i=1,2,\cdots,n~$ any two points of the lattice in
 $\mathbb{R}^{n-i+1}$ with basis $$(A_i,0,...,0),
 (a_{i+1,i},A_{i+1},0,\cdots,0), \cdots,(a_{n,i},a_{n,i+1},
\cdots,a_{n,n-1},A_n)$$ are at a distance at least $A_i$ apart. Here we shall be considering the following conjecture of Woods :\\

\noindent {\bf Conjecture (Woods).} If $A_1A_2\cdots
A_n = 1$ and $A_i\leq A_1$ for each $i$ then any closed sphere in
$\mathbb{R}^n$ of radius $\sqrt{n}/2$ contains a point of $\mathbb{L}.$\vspace{2mm}

\noindent This conjecture is known to be true for $n \leq 9$. Woods \cite{W1, W2, W3} proved it for $4\leq n\leq 6~$.  Hans-Gill et al. \cite{HRS7, HRS8} proved it for $n=7$ and $n=8$. In a previous paper, the authors \cite{KR1} proved it for $n=9$.
In \cite{KR2}, the authors  have  obtained estimates to the Conjecture
of  Woods for $10 \le n \le 33$. In particular we obtained a weaker result for $n=10$
that if hypothesis of Woods' Conjecture holds, then any closed sphere in $\mathbb{R}^{10}$ of radius
$\frac{\sqrt{10.3}}{2}$
 contains a point of $\mathbb{L}.$ In 2017,  Regev et al. \cite{RSW} showed that Woods' Conjecture is not always true, it is false
   for $ n\ge 30$. It will be of great interest to find  the largest (smallest)
value of $n$, $10 \le n \le 29$, for which Woods' Conjecture is true (false). In this direction, we find that Woods' Conjecture is true for $n = 10.$ \vspace{2mm}

\noindent Woods \cite{W2,W3} showed that his conjecture implies the following  conjecture:\vspace{3mm}

\noindent {\bf Conjecture I.} If $\wedge$ is a lattice of determinant 1 and there is a sphere $|X|<R$ which contains no point of $\wedge$ other than $O$ and has $n$ linearly independent points of $\wedge$ on its boundary then $\wedge$ is a covering lattice for the closed sphere of radius $ \sqrt{n/4}.$ Equivalently every closed sphere of radius $\sqrt{n/4}$ lying in $\mathbb{R}^n$ contains a point of $\wedge$.\vspace{2mm}

 It is well known that together with the result of McMullen \cite{Mc}, truth of  Conjecture I for a fixed $n$ would imply the following long standing classical conjecture attributed to Minkowski  on the product of $n$ non-homogeneous linear forms in $n$ variables:\\

{\noindent \bf  Conjecture (Minkowski).} Let $L_i = a_{i1}x_1+\cdots + a_{in}x_n,$ $1\leq i\leq n,$ be
$n$ real linear forms in $n$ variables $x_1,\cdots,x_n$  having determinant $\Delta = \det{(a_{ij})}\neq 0.$ For any
given real numbers $c_1,\cdots,c_n$ there exist integers $x_1,\cdots,x_n$ such that
\begin{equation*}|(L_1+c_1)\cdots(L_n+c_n)|\leq |\Delta|/2^n.\end{equation*}
Minkowski's Conjecture is known to be true for $n \leq 9.$  For more detailed history of  {\mm}'s Conjecture and related results, see Gruber \cite{PG}, Gruber and Lekkerkerker
\cite{GL}, Bambah et al. \cite{BDH} and Hans-Gill et al. \cite{HRS7}. While answering a question of Shapira and Weiss \cite{SW}, the
authors along with Hans-Gill \cite{KHR} have given another proof of Minkowski's Conjecture  for $n \le 7$.\vspace{1mm}

\noindent In this paper we shall prove

\begin{theorem} Woods' Conjecture is true for n = 10.\end{theorem}

Therefore Conjecture I and hence   Minkowski's Conjecture is proved for $n=10$. One notes that falsehood of Woods' Conjecture for $n \ge 30$ does not mean that Minkowski's Conjecture is also false for those $n$.\vspace{2mm}

 We use the notations and method of proof of Hans-Gill et al. \cite{HRS7, HRS8}. In this method, we need to maximize /minimize
frequently functions of several variables. For n=7 and 8, Hans-Gill et al. \cite{HRS7, HRS8} did all the calculations by hand using calculus only. While proving it for $n=9$, see \cite{KR1}, we reduced the number of variables  one by one using calculus and replaced them with values where it
 could have its optimum value. This way we obtained functions
in at most 3 variables. Finally we  arrived at a conclusion by plotting 2 or 3 dimensional
graphs in software Mathematica. Because of the heavy calculation work, the
proof became very lengthy. The detailed proof of Woods' Conjecture for $n=9$ consisting of 132 pages can be seen at
arXiv:1410.5743v1[math.NT]. With the increase of one more dimension i.e. for $n=10,$ the same process is extremely difficult to give a result without using some  computational package. Here we use  Optimization tools of the software Mathematica (non-linear global optimization) to prove Woods' Conjecture for $n=10$.

All these optimization computations were initially verified by package Lingo. The authors are very grateful to Lindo Systems for providing its  access free of charge.

\section{Preliminary Lemmas}
For a unit sphere $S_{n}$ with center $O$ in $\mathbb{R}^{n}$, let $\Delta(S_{n})$ be the \emph{critical determinant} of $S_{n}$, defined as
\vspace{-2mm}$$\Delta(S_{n}) = \inf\{d(\Lambda):\Lambda~\mbox{has no non-zero point in the interior of}~S_{n}\},$$ where $d(\Lambda)$ denotes the determinant of the lattice $\Lambda$.
 \par Let $\mathbb{L}$ be a lattice in $\mathbb{R}^{n}$ reduced in the sense of Korkine and Zolotareff and $A_1, A_2,\ldots,A_n$ be defined as in Section 1. We state below some preliminary lemmas. Lemmas 1 and 2 are due to Woods \cite{W1}, Lemma 3 is due to Korkine and Zolotareff \cite{KZ} and Lemma 4 is due to  Pendavingh and Van Zwam \cite{PV}.
 In Lemma 5, the cases $n=2$ and $3$ are classical results of Lagrange and
 Gauss; $n=4$ and $5$ are due to Korkine and Zolotareff \cite{KZ} while $n=6, 7$ and $8$ are due to Blichfeldt \cite{Bh}.\vspace{2mm}\\
\begin{lemma} \emph{If $2\Delta(S_{n+1})A_{1}^{n}\geq d(\mathbb{L})$, then any closed sphere of radius
$$R=A_{1}\{ 1-(A_{1}^{n} \Delta (S_{n+1})/ d(\mathbb{L}))^{2}\}^{1/2}$$ in $\mathbb{R}^{n}$ contains a point of $\mathbb{L}$.}\end{lemma}
\begin{lemma}   \emph{For a fixed  integer $i$ with $1\leqslant i\leqslant n-1$, denote by $\mathbb{L}_{1}$ the lattice in $\mathbb{R}^{i}$ with the reduced basis $$(A_1,0,0,\ldots,0),(a_{2,1},A_2,0,\ldots,0),\ldots,(a_{i,1},a_{i,2},\ldots,a_{i,i-1},A_i) $$ and denote by $\mathbb{L}_{2}$ the lattice in $\mathbb{R}^{n-i}$ with the reduced basis $$(A_{i+1},0,0,\ldots,0),(a_{i+2,i+1},A_{i+2},0,\ldots,0),\ldots,(a_{n,i+1},a_{n,i+2},\ldots,a_{n,n-1},A_n). $$
If any sphere in $\mathbb{R}^{i}$ of radius $r_{1}$ contains a point of $\mathbb{L}_{1}$ and if any sphere in $\mathbb{R}^{n-i}$ of radius $r_{2}$ contains a point
of $\mathbb{L}_{2}$ then any sphere in $\mathbb{R}^{n}$ of radius $(r_{1}^{2}+r_{2}^{2})^{1/2}$ contains a point of $\mathbb{L}$.}\end{lemma}
\begin{lemma}  \emph{ For all relevant $i$, $A_{i+1}^2\geq\frac{3}{4}A_i^2$ and $A_{i+2}^2\geq\frac{2}{3}A_i^2
~$.}\end{lemma}
\begin{lemma}  \emph{ For all relevant $i$, $A_{i+4}^2\geq0.46873A_i^2 $~.}\end{lemma}
\begin{lemma}  \emph{$~~~\Delta (S_n) = ~1/\sqrt{2}, ~1/2, ~1/2\sqrt{2}, ~\sqrt{3}/8,1/8$ and $1/16  ~for~
 n=~3,~4,~5,$ $6,~7$ and $~8$ ~respectively.}\end{lemma}
\section{Plan of the Proof}\label{Plan} As remarked earlier, we use the notation and approach of Hans-Gill et al. \cite{HRS7} and that of \cite{HRS8}. We include some of the details given there for the convenience of the reader. We assume that Woods' Conjecture is false for $n=10$ and derive a contradiction. Let $\mathbb{L}$ be lattice satisfying the hypothesis of the conjecture for $n=10$ i.e.  $A_{1}A_{2}\cdots A_{10}=1$~\mbox{and}~ $A_{i}\leqslant A_{1}$ for each $i$. Suppose that there exists a closed sphere of radius $\sqrt{10}/2$ in $\mathbb{R}^{10}$ that contains no point of $\mathbb{L}$. Write $A=A_1^2,~ B=A_2^2,~ C=A_3^2,\ldots, J=A_{10}^2$. So we have $ABCDEFGHIJ=1$.\vspace{2mm}\\
$~~~$ If
$(\lambda_{1},\lambda_{2},\cdots,\lambda_{s})$ is an ordered partition of $n$, then the conditional inequality arising
from it, by using Lemmas 1 and 2, is also denoted by $(\lambda_{1},\lambda_{2},\cdots,\lambda_{s})$. If the conditions
in an inequality $(\lambda_{1},\lambda_{2},\cdots,\lambda_{s})$ are satisfied then we say that
$(\lambda_{1},\lambda_{2},\cdots,\lambda_{s})$ holds.\vspace{2mm}\\
For example the inequality (1, 1, 1, 1, 1, 1, 1, 1, 2) results in the conditional inequality :
\begin{equation}\label{2}
{\rm if~~} 2I\geq J~~ {\rm then~~~} A+B+C+D+E+F+G+H+4I-\frac{2I^2}{J}>10.
\end{equation}
Since $4I-2I^2/J\leq 2J$,  the second inequality in (\ref{2}) gives
\begin{equation}\label{3}
A+B+C+D+E+F+G+H+2J>10.
\end{equation}
One may remark here that the condition $2I\geq J$ is necessary only if we want to use inequality (\ref{2}), but it is not necessary if we want to  use the weaker inequality (\ref{3}). This is so because if $2I< J$, using the partition $(1,1)$ in place of $(2)$ for the relevant part, we get the upper bound $I+J$ which is clearly less than $2J$. We shall call inequalities of type (\ref{3}) as \emph{weak} inequalities and inequalities of type (\ref{2}) as \emph{strong} inequalities.\vspace{1mm} \\
$~~~$Sometimes, instead of Lemma 1, we are able to use the fact that Woods Conjecture is true for dimensions less than or equal to $9$. The
use of this is indicated by putting $^{*}$ on the corresponding part of the partition. For example, the inequality $(6^{*},4)$
 is  \begin{equation} {\rm if~~} G^4ABCDEF\geq 2~~ {\rm then~~~} 6(ABCDEF)^{\frac{1}{6}}+4G-\frac{1}{2}G^{5}ABCDEF > 10,\end{equation}
the hypothesis of the conjecture in $6$ variables being satisfied.\vspace{2mm}\\
$~~~~$We observe that the inequalities of the type $(3,1,1,\cdots,1)$, $(2,1,1,\cdots,1)$, $(1,2,1,\cdots,1)$, $(1,2,2,1,\cdots,1)$, $(3,2,1,\cdots,1)$ etc. always hold. \vspace{2mm}\\

$~~~~~$We can assume $A>1$, because if $A\leq 1$, we must
have  $A=B=C=D=E=F=G=H=I=J=1.$ In this case Woods' Conjecture  can be seen to be true using inequality $(1,1,1,1,1,1,1,1,1,1)$. Also it is known that $A\leq \gamma_{10}<2.2636302$ (See \cite{CE, CS}), where $\gamma_n$ denotes the Hermite's constant for minima of positive definite quadratic forms.\vspace{2mm}\\
$~~~~~$Each of $B,C,\ldots,J$ can either be $>1$ or $\leqslant$ 1. This give rise to $2^9=512$ cases. Case 1, where each of $B,C,D\cdots,J > 1$ does not arise as $ABCDEFGHIJ=1$. In Section \ref{Easy}, the 353 easy cases have been considered under Propositions 1-4. The remaining 158 cases need much more intricate analysis of the available inequalities. Out of these 158 cases, 151  cases are somewhat less difficult and have been dealt in Section \ref{Difficult}. The remaining 7 cases are very difficult to solve as the ranges of the variables have to be divided into many sub-intervals. So these cases have been dealt seperately in Section \ref{Most Difficult}. For the cases in Sections \ref{Difficult} and  \ref{Most Difficult}, we have used the Optimization tool of the software Mathematica(non-linear global optimization) for the optimization of the functions in 10 variables $A, B, C, D, E, F, G, H, I, J$ under a number of constraints.\vspace{1mm}\\
$~~~~~$We write all the 512 cases in lexicographical order. For the 353 easy cases, the inequalities used to get a contradiction and the propositions where they are dealt with are listed in Table I. For 158 difficult cases either the proposition where these are dealt with  or the main  inequalities used to get a contradiction are  also listed in Table I.\\
\textbf{Remark:} In many cases there are alternative ways to
get a contradiction. We have chosen to describe the method which we find
convenient.
{\footnotesize$${\rm \bf Table~~~I}$$
\begin{tabular}{lcccccccccccc}
Case & A & B & C & D & E & F & G & H & I &J&Proposition& Inequalities \\ &&&&&&&&&&&&\\
1&$>$&$>$&$>$&$>$&$>$&$>$&$>$&$>$&$>$&$>$&$-$&ABCDEFGHIJ=1\\
2&$>$&$>$&$>$&$>$&$>$&$>$&$>$&$>$&$>$&$\leq$&$1$&$(8^*,2)$\\
3&$>$&$>$&$>$&$>$&$>$&$>$&$>$&$>$&$\leq$ &$>$&3(i)&$(1,1,1,1,1,1,1,2,1)$ \\
4&$>$&$>$&$>$&$>$&$>$&$>$&$>$&$>$&$\leq$ &$\leq$&1&$(7^*,3)$ \\
5&$>$&$>$&$>$&$>$&$>$&$>$&$>$&$\leq$&$>$ &$>$&3(i)&$(1,1,1,1,1,1,2,1,1)$\\
6&$>$&$>$&$>$&$>$&$>$&$>$&$>$&$\leq$&$>$ &$\leq$&1&$(8^*,2)$\\
7&$>$&$>$&$>$&$>$&$>$&$>$&$>$&$\leq$&$\leq$ &$>$&3(v)&$(1,1,1,1,1,1,3,1)$ \\
8&$>$&$>$&$>$&$>$&$>$&$>$&$>$&$\leq$&$\leq$ &$\leq$&6&$-$ \\
9&$>$&$>$&$>$&$>$&$>$&$>$&$\leq$&$>$&$>$&$>$&3(i)&$(1,1,1,1,1,2,1,1,1)$ \\
10&$>$&$>$&$>$&$>$&$>$&$>$&$\leq$&$>$&$>$&$\leq$&1&$(8^*,2)$ \\
11&$>$&$>$&$>$&$>$&$>$&$>$&$\leq$&$>$&$\leq$ &$>$&3(ii)&$(1,1,1,1,2,2,1)$ \\
12&$>$&$>$&$>$&$>$&$>$&$>$&$\leq$&$>$&$\leq$ &$\leq$&1&$(7^*,3)$ \\
13&$>$&$>$&$>$&$>$&$>$&$>$&$\leq$&$\leq$&$>$ &$>$&3(v)&$(1,1,1,1,1,3,1,1)$ \\
14&$>$&$>$&$>$&$>$&$>$&$>$&$\leq$&$\leq$&$>$ &$\leq$&1&$(8^*,2)$ \\
15&$>$&$>$&$>$&$>$&$>$&$>$&$\leq$&$\leq$&$\leq$ &$>$&$6$ &$-$ \\
16&$>$&$>$&$>$&$>$&$>$&$>$&$\leq$&$\leq$&$\leq$ &$\leq$&$7$ &$-$ \\
17&$>$&$>$&$>$&$>$&$>$&$\leq$&$>$&$>$&$>$ &$>$&3(i)&$(1,1,1,1,2,1,1,1,1)$ \\
18&$>$&$>$&$>$&$>$&$>$&$\leq$&$>$&$>$&$>$ &$\leq$&1&$(8^*,2)$ \\
19&$>$&$>$&$>$&$>$&$>$&$\leq$&$>$&$>$&$\leq$ &$>$&3(ii)&$(1,1,1,1,2,1,2,1)$ \\
20&$>$&$>$&$>$&$>$&$>$&$\leq$&$>$&$>$&$\leq$ &$\leq$&1&$(7^*,3)$ \\
21&$>$&$>$&$>$&$>$&$>$&$\leq$&$>$&$\leq$&$>$ &$>$&3(ii)&$(1,1,1,1,2,2,1,1)$ \\
22&$>$&$>$&$>$&$>$&$>$&$\leq$&$>$&$\leq$&$>$ &$\leq$&1&$(8^*,2)$ \\
23&$>$&$>$&$>$&$>$&$>$&$\leq$&$>$&$\leq$&$\leq$ &$>$&3(ix)& $(1,1,1,1,2,3,1)$ \\
24&$>$&$>$&$>$&$>$&$>$&$\leq$&$>$&$\leq$&$\leq$ &$\leq$&-& $(2,2,2,3,1)$ \\
25&$>$&$>$&$>$&$>$&$>$&$\leq$&$\leq$&$>$&$>$ &$>$&3(v)& $(1,1,1,1,3,1,1,1)$\\
26&$>$&$>$&$>$&$>$&$>$&$\leq$&$\leq$&$>$&$>$ &$\leq$&1& $(8^*,2)$\\
27&$>$&$>$&$>$&$>$&$>$&$\leq$&$\leq$&$>$&$\leq$ &$>$&3(xi)& $(3,1,3,2,1)$ \\
28&$>$&$>$&$>$&$>$&$>$&$\leq$&$\leq$&$>$&$\leq$ &$\leq$&1& $(7^*,3)$ \\
29&$>$&$>$&$>$&$>$&$>$&$\leq$&$\leq$&$\leq$&$>$ &$>$&6&$ -$\\
30&$>$&$>$&$>$&$>$&$>$&$\leq$&$\leq$&$\leq$&$>$ &$\leq$&1&$ (8^*,2)$\\
31&$>$&$>$&$>$&$>$&$>$&$\leq$&$\leq$&$\leq$&$\leq$ &$>$&7& $-$\\
32&$>$&$>$&$>$&$>$&$>$&$\leq$&$\leq$&$\leq$&$\leq$ &$\leq$&7& $-$\\
33&$>$&$>$&$>$&$>$&$\leq$&$>$&$>$&$>$&$>$ &$>$&3(i)&$(1,1,1,2,1,1,1,1,1)$ \\
34&$>$&$>$&$>$&$>$&$\leq$&$>$&$>$&$>$&$>$ &$\leq$&1&$(8^*,2)$ \\
35&$>$&$>$&$>$&$>$&$\leq$&$>$&$>$&$>$&$\leq$ &$>$&3(ii)& $(1,1,1,2,1,1,2,1)$ \\
36&$>$&$>$&$>$&$>$&$\leq$&$>$&$>$&$>$&$\leq$ &$\leq$&1& $(7^*,3)$ \\
37&$>$&$>$&$>$&$>$&$\leq$&$>$&$>$&$\leq$&$>$ &$>$&3(ii)&$(1,1,1,2,1,2,1,1)$ \\
38&$>$&$>$&$>$&$>$&$\leq$&$>$&$>$&$\leq$&$>$ &$\leq$&1&$(8^*,2)$ \\
39&$>$&$>$&$>$&$>$&$\leq$&$>$&$>$&$\leq$&$\leq$ &$>$&3(xi)&$(3,2,1,3,1)$ \\
40&$>$&$>$&$>$&$>$&$\leq$&$>$&$>$&$\leq$&$\leq$ &$\leq$&-&$(1,2,2,1,3,1)$ \\
\end{tabular}
\newpage
 \begin{tabular}{ccccccccccccc}
 Case & A & B & C & D & E & F & G & H & I &J&Proposition& Inequalities \\ &&&&&&&&&&&&\\
41&$>$&$>$&$>$&$>$&$\leq$&$>$&$\leq$&$>$&$>$ &$>$&3(ii)&$(1,1,1,2,2,1,1,1)$\\
42&$>$&$>$&$>$&$>$&$\leq$&$>$&$\leq$&$>$&$>$ &$\leq$&1&$(8^*,2)$\\
43&$>$&$>$&$>$&$>$&$\leq$&$>$&$\leq$&$>$&$\leq$ &$>$&3(iii)& $(1,1,1,2,2,2,1)$ \\
44&$>$&$>$&$>$&$>$&$\leq$&$>$&$\leq$&$>$&$\leq$ &$\leq$&1& $(7^*,3)$ \\
45&$>$&$>$&$>$&$>$&$\leq$&$>$&$\leq$&$\leq$&$>$ &$>$&3(xi)&$(3,2,3,1,1)$ \\
46&$>$&$>$&$>$&$>$&$\leq$&$>$&$\leq$&$\leq$&$>$ &$\leq$&1&$(8^*,2)$ \\
47&$>$&$>$&$>$&$>$&$\leq$&$>$&$\leq$&$\leq$&$\leq$ &$>$&-& $(1,2,2,3,1,1)$\\
48&$>$&$>$&$>$&$>$&$\leq$&$>$&$\leq$&$\leq$&$\leq$ &$\leq$&6& $-$\\
49&$>$&$>$&$>$&$>$&$\leq$&$\leq$&$>$&$>$&$>$ &$>$&3(v)&$(1,1,1,3,1,1,1,1)$ \\
50&$>$&$>$&$>$&$>$&$\leq$&$\leq$&$>$&$>$&$>$ &$\leq$&1&$(8^*,2)$ \\
51&$>$&$>$&$>$&$>$&$\leq$&$\leq$&$>$&$>$&$\leq$ &$>$&3(xi)&$(3,3,1,2,1)$\\
52&$>$&$>$&$>$&$>$&$\leq$&$\leq$&$>$&$>$&$\leq$ &$\leq$&1&$(7^*,3)$\\
53&$>$&$>$&$>$&$>$&$\leq$&$\leq$&$>$&$\leq$&$>$ &$>$&3(xi)&$(3,3,2,1,1)$\\
54&$>$&$>$&$>$&$>$&$\leq$&$\leq$&$>$&$\leq$&$>$ &$\leq$&1&$(8^*,2)$\\
55&$>$&$>$&$>$&$>$&$\leq$&$\leq$&$>$&$\leq$&$\leq$&$>$&3(vii)&$(3,3,3,1)$\\
56&$>$&$>$&$>$&$>$&$\leq$&$\leq$&$>$&$\leq$&$\leq$&$\leq$&-&$(1,2,3,3,1)$\\
57&$>$&$>$&$>$&$>$&$\leq$&$\leq$&$\leq$&$>$&$>$ &$>$&-&$(1,2,3,1,2,1)$\\
58&$>$&$>$&$>$&$>$&$\leq$&$\leq$&$\leq$&$>$&$>$ &$\leq$&1&$(8^*,2)$\\
59&$>$&$>$&$>$&$>$&$\leq$&$\leq$&$\leq$&$>$&$\leq$ &$>$&-&$(1,2,3,1,2,1)$ \\
60&$>$&$>$&$>$&$>$&$\leq$&$\leq$&$\leq$&$>$&$\leq$ &$\leq$&1&$(7^*,3)$ \\
61&$>$&$>$&$>$&$>$&$\leq$&$\leq$&$\leq$&$\leq$&$>$ &$>$& 7&$-$\\
62&$>$&$>$&$>$&$>$&$\leq$&$\leq$&$\leq$&$\leq$&$>$ &$\leq$&1&$(8^*,2)$\\
63&$>$&$>$&$>$&$>$&$\leq$&$\leq$&$\leq$&$\leq$&$\leq$&$>$&7&$-$\\
64&$>$&$>$&$>$&$>$&$\leq$&$\leq$&$\leq$&$\leq$&$\leq$&$\leq$&7&$-$\\
65 &$>$&$>$&$>$&$\leq$&$>$&$>$&$>$&$>$&$>$&$>$&3(i)&$(1,1,2,1,1,1,1,1,1)$\\
66 &$>$&$>$&$>$&$\leq$&$>$&$>$&$>$&$>$&$>$&$\leq$&1&$(8^*,2)$\\
67&$>$&$>$&$>$&$\leq$&$>$&$>$&$>$&$>$&$\leq$&$>$&3(ii)&$(1,1,2,1,1,1,2,1)$ \\
68&$>$&$>$&$>$&$\leq$&$>$&$>$&$>$&$>$&$\leq$&$\leq$&1&$(7^*,3)$ \\
69&$>$&$>$&$>$&$\leq$&$>$&$>$&$>$&$\leq$&$>$&$>$&3(ii)&$(1,1,2,1,1,2,1,1)$\\
70&$>$&$>$&$>$&$\leq$&$>$&$>$&$>$&$\leq$&$>$&$\leq$&1&$(8^*,2)$\\
71&$>$&$>$&$>$&$\leq$&$>$&$>$&$>$&$\leq$&$\leq$&$>$&3(ix)&$(2,2,1,1,3,1)$\\
72&$>$&$>$&$>$&$\leq$&$>$&$>$&$>$&$\leq$&$\leq$&$\leq$&-&$(2,2,1,1,3,1)$\\
73&$>$&$>$&$>$&$\leq$&$>$&$>$&$\leq$&$>$&$>$&$>$&3(ii)&$(1,1,2,1,2,1,1,1)$\\
74&$>$&$>$&$>$&$\leq$&$>$&$>$&$\leq$&$>$&$>$&$\leq$&1&$(8^*,2)$\\
75&$>$&$>$&$>$&$\leq$&$>$&$>$&$\leq$&$>$&$\leq$&$>$&3(iii)&$(1,1,2,1,2,2,1)$\\
76&$>$&$>$&$>$&$\leq$&$>$&$>$&$\leq$&$>$&$\leq$&$\leq$&1&$(7^*,3)$\\
77&$>$&$>$&$>$&$\leq$&$>$&$>$&$\leq$&$\leq$&$>$&$>$&3(ix)&$(2,2,1,3,1,1)$\\
78&$>$&$>$&$>$&$\leq$&$>$&$>$&$\leq$&$\leq$&$>$&$\leq$&1&$(8^*,2)$\\
79&$>$&$>$&$>$&$\leq$&$>$&$>$&$\leq$&$\leq$&$\leq$&$>$&-&$(2,2,1,3,1,1)$\\
80&$>$&$>$&$>$&$\leq$&$>$&$>$&$\leq$&$\leq$&$\leq$&$\leq$&-&$(2,2,1,3,1,1)$\\
\end{tabular}
\newpage
 \begin{tabular}{ccccccccccccc}
 Case & A & B & C & D & E & F & G & H & I &J&Proposition& Inequalities \\ &&&&&&&&&&&&\\
81&$>$&$>$&$>$&$\leq$&$>$&$\leq$&$>$&$>$&$>$&$>$&3(ii)&$(1,1,2,2,1,1,1,1)$\\
82&$>$&$>$&$>$&$\leq$&$>$&$\leq$&$>$&$>$&$>$&$\leq$&1&$(8^*,2)$\\
83&$>$&$>$&$>$&$\leq$&$>$&$\leq$&$>$&$>$&$\leq$&$>$&3(iii)&$(1,1,2,2,1,2,1)$\\
84&$>$&$>$&$>$&$\leq$&$>$&$\leq$&$>$&$>$&$\leq$&$\leq$&1&$(7^*,3)$\\
85&$>$ &$>$&$>$&$\leq$&$>$&$\leq$&$>$&$\leq$&$>$&$>$&3(iii)&$(1,1,2,2,2,1,1)$\\
86&$>$ &$>$&$>$&$\leq$&$>$&$\leq$&$>$&$\leq$&$>$&$\leq$&1&$(8^*,2)$\\
87&$>$&$>$&$>$&$\leq$&$>$&$\leq$&$>$&$\leq$&$\leq$&$>$&3(x)&$(2,2,2,3,1)$\\
88&$>$&$>$&$>$&$\leq$&$>$&$\leq$&$>$&$\leq$&$\leq$&$\leq$&-&$(2,2,2,3,1)$\\
89&$>$&$>$&$>$&$\leq$&$>$&$\leq$&$\leq$&$>$&$>$&$>$&3(ix)&$(2,2,3,1,1,1)$\\
90&$>$&$>$&$>$&$\leq$&$>$&$\leq$&$\leq$&$>$&$>$&$\leq$&1&$(8^*,2)$\\
91&$>$&$>$&$>$&$\leq$&$>$&$\leq$&$\leq$&$>$&$\leq$&$>$&3(x)&$(2,2,3,2,1)$\\
92&$>$&$>$&$>$&$\leq$&$>$&$\leq$&$\leq$&$>$&$\leq$&$\leq$&1&$(7^*,3)$\\
93&$>$&$>$&$>$&$\leq$&$>$&$\leq$&$\leq$&$\leq$&$>$&$>$&-& $(2,2,3,1,1,1)$\\
94&$>$&$>$&$>$&$\leq$&$>$&$\leq$&$\leq$&$\leq$&$>$&$\leq$&1& $(8^*,2)$\\
95&$>$&$>$&$>$&$\leq$&$>$&$\leq$&$\leq$&$\leq$&$\leq$&$>$&6&$-$\\
96&$>$&$>$&$>$&$\leq$&$>$&$\leq$&$\leq$&$\leq$&$\leq$&$\leq$&-&$(2,2,3,1,1,1)$\\
97&$>$&$>$&$>$&$\leq$&$\leq$&$>$&$>$&$>$&$>$&$>$&3(v)&$(1,1,3,1,1,1,1,1)$\\
98&$>$&$>$&$>$&$\leq$&$\leq$&$>$&$>$&$>$&$>$&$\leq$&1&$(8^*,2)$\\
99&$>$&$>$&$>$&$\leq$&$\leq$&$>$&$>$&$>$&$\leq$&$>$&3(ix)&$(2,3,1,1,2,1)$\\
100&$>$&$>$&$>$&$\leq$&$\leq$&$>$&$>$&$>$&$\leq$&$\leq$&1&$(7^*,3)$\\
101&$>$&$>$&$>$&$\leq$&$\leq$&$>$&$>$&$\leq$&$>$&$>$&3(ix)&$(2,3,1,2,1,1)$\\
102&$>$&$>$&$>$&$\leq$&$\leq$&$>$&$>$&$\leq$&$>$&$\leq$&1&$(8^*,2)$\\
103&$>$ &$>$&$>$&$\leq$&$\leq$&$>$&$>$&$\leq$&$\leq$&$>$&3(xi)&$(2,3,1,3,1)$ \\
104&$>$ &$>$&$>$&$\leq$&$\leq$&$>$&$>$&$\leq$&$\leq$&$\leq$&-&$(2,3,1,3,1)$ \\
105&$>$&$>$&$>$&$\leq$&$\leq$&$>$&$\leq$&$>$&$>$&$>$&3(ix)& $(2,3,2,1,1,1)$\\
106&$>$&$>$&$>$&$\leq$&$\leq$&$>$&$\leq$&$>$&$>$&$\leq$&1& $(8^*,2)$\\
107&$>$&$>$&$>$&$\leq$&$\leq$&$>$&$\leq$&$>$&$\leq$&$>$&3(x)&$(2,3,2,2,1)$\\
108&$>$&$>$&$>$&$\leq$&$\leq$&$>$&$\leq$&$>$&$\leq$&$\leq$&1&$(7^*,3)$\\
109&$>$&$>$&$>$&$\leq$&$\leq$&$>$&$\leq$&$\leq$&$>$&$>$&3(xi)&$(2,3,3,1,1)$\\
110&$>$&$>$&$>$&$\leq$&$\leq$&$>$&$\leq$&$\leq$&$>$&$\leq$&1&$(8^*,2)$\\
111&$>$&$>$&$>$&$\leq$&$\leq$&$>$&$\leq$&$\leq$&$\leq$&$>$&-&$(2,2,1,2,1,1,1)$\\
112&$>$&$>$&$>$&$\leq$&$\leq$&$>$&$\leq$&$\leq$&$\leq$&$\leq$&-&$(2,2,1,2,1,1,1)$\\
113&$>$&$>$&$>$&$\leq$&$\leq$&$\leq$&$>$&$>$&$>$&$>$&-&$(2,2,2,1,1,1,1)$\\
114&$>$&$>$&$>$&$\leq$&$\leq$&$\leq$&$>$&$>$&$>$&$\leq$&1&$(8^*,2)$\\
115&$>$&$>$&$>$&$\leq$&$\leq$&$\leq$&$>$&$>$&$\leq$&$>$&-&$(2,2,2,1,2,1)$\\
116&$>$&$>$&$>$&$\leq$&$\leq$&$\leq$&$>$&$>$&$\leq$&$\leq$&1&$(7^*,3)$\\
117&$>$ &$>$&$>$&$\leq$&$\leq$&$\leq$&$>$&$\leq$&$>$&$>$&-&$(2,2,2,2,1,1)$\\
118&$>$ &$>$&$>$&$\leq$&$\leq$&$\leq$&$>$&$\leq$&$>$&$\leq$&1&$(8^*,2)$\\
119&$>$&$>$&$>$&$\leq$&$\leq$&$\leq$&$>$&$\leq$&$\leq$&$>$&-&$(2,2,2,2,1,1)$\\
120&$>$&$>$&$>$&$\leq$&$\leq$&$\leq$&$>$&$\leq$&$\leq$&$\leq$&-&$(2,2,2,2,1,1)$\\
\end{tabular}
\newpage
 \begin{tabular}{ccccccccccccc}
 Case & A & B & C & D & E & F & G & H & I &J&Proposition& Inequalities \\ &&&&&&&&&&&&\\
121&$>$&$>$&$>$&$\leq$&$\leq$&$\leq$&$\leq$&$>$&$>$&$>$&6&$-$\\
122&$>$&$>$&$>$&$\leq$&$\leq$&$\leq$&$\leq$&$>$&$>$&$\leq$&1&$(8^*,2)$\\
123&$>$&$>$&$>$&$\leq$&$\leq$&$\leq$&$\leq$&$>$&$\leq$&$>$&6&$-$\\
124&$>$&$>$&$>$&$\leq$&$\leq$&$\leq$&$\leq$&$>$&$\leq$&$\leq$&1&$(7^*,3)$\\
125&$>$&$>$&$>$&$\leq$&$\leq$&$\leq$&$\leq$&$\leq$&$>$ &$>$&7&$-$\\
126&$>$&$>$&$>$&$\leq$&$\leq$&$\leq$&$\leq$&$\leq$&$>$ &$\leq$&1&$(8^*,2)$\\
127&$>$&$>$&$>$&$\leq$&$\leq$&$\leq$&$\leq$&$\leq$&$\leq$&$>$&6&$-$\\
128&$>$&$>$&$>$&$\leq$&$\leq$&$\leq$&$\leq$&$\leq$&$\leq$&$\leq$&-&$(1,4,1,\cdots,1),(1,2,1,\cdots,1)$\\
129&$>$&$>$&$\leq$&$>$&$>$&$>$&$>$&$>$&$>$&$>$&3(i)&$(1,2,1,1,1,1,1,1,1)$ \\
130&$>$&$>$&$\leq$&$>$&$>$&$>$&$>$&$>$&$>$&$\leq$&1&$(8^*,2)$ \\
131&$>$&$>$&$\leq$&$>$&$>$&$>$&$>$&$>$&$\leq$&$>$&3(ii)&$(1,2,1,1,1,1,2,1)$ \\
132&$>$&$>$&$\leq$&$>$&$>$&$>$&$>$&$>$&$\leq$&$\leq$&1&$(7^*,3)$ \\
133&$>$&$>$&$\leq$&$>$&$>$&$>$&$>$&$\leq$&$>$&$>$&3(ii)&$(1,2,1,1,1,2,1,1)$ \\
134&$>$&$>$&$\leq$&$>$&$>$&$>$&$>$&$\leq$&$>$&$\leq$&1&$(8^*,2)$ \\
135&$>$&$>$&$\leq$&$>$&$>$&$>$&$>$&$\leq$&$\leq$&$>$&3(vi)&$(3,1,1,1,3,1)$ \\
136&$>$&$>$&$\leq$&$>$&$>$&$>$&$>$&$\leq$&$\leq$&$\leq$&-&$(1,2,1,1,1,3,1)$ \\
137&$>$&$>$&$\leq$&$>$&$>$&$>$&$\leq$&$>$&$>$&$>$&3(ii)&$(1,2,1,1,2,1,1,1)$ \\
138&$>$&$>$&$\leq$&$>$&$>$&$>$&$\leq$&$>$&$>$&$\leq$&1&$(8^*,2)$ \\
139&$>$&$>$&$\leq$&$>$&$>$&$>$&$\leq$&$>$&$\leq$&$>$&3(iii)&$(1,2,1,1,2,2,1)$ \\
140&$>$&$>$&$\leq$&$>$&$>$&$>$&$\leq$&$>$&$\leq$&$\leq$&1&$(7^*,3)$ \\
141&$>$&$>$&$\leq$&$>$&$>$&$>$&$\leq$&$\leq$&$>$&$>$&3(vi)&$(3,1,1,3,1,1)$ \\
142&$>$&$>$&$\leq$&$>$&$>$&$>$&$\leq$&$\leq$&$>$&$\leq$&1&$(8^*,2)$ \\
143&$>$&$>$&$\leq$&$>$&$>$&$>$&$\leq$&$\leq$&$\leq$&$>$&-&$(1,2,1,1,3,1,1)$ \\
144&$>$&$>$&$\leq$&$>$&$>$&$>$&$\leq$&$\leq$&$\leq$&$\leq$&-&$(1,2,1,1,3,1,1)$ \\
145&$>$&$>$&$\leq$&$>$&$>$&$\leq$&$>$&$>$&$>$&$>$&3(ii)&$(1,2,1,2,1,1,1,1)$ \\
146&$>$&$>$&$\leq$&$>$&$>$&$\leq$&$>$&$>$&$>$&$\leq$&1&$(8^*,2)$ \\
147&$>$&$>$&$\leq$&$>$&$>$&$\leq$&$>$&$>$&$\leq$&$>$&3(iii)&$(1,2,1,2,1,2,1)$ \\
148&$>$&$>$&$\leq$&$>$&$>$&$\leq$&$>$&$>$&$\leq$&$\leq$&1&$(7^*,3)$ \\
149&$>$&$>$&$\leq$&$>$&$>$&$\leq$&$>$&$\leq$&$>$&$>$&3(iii)&$(1,2,1,2,2,1,1)$ \\
150&$>$&$>$&$\leq$&$>$&$>$&$\leq$&$>$&$\leq$&$>$&$\leq$&&$(8^*,2)$ \\
151&$>$&$>$&$\leq$&$>$&$>$&$\leq$&$>$&$\leq$&$\leq$&$>$&3(xi)&$(3,1,2,3,1)$\\
152&$>$&$>$&$\leq$&$>$&$>$&$\leq$&$>$&$\leq$&$\leq$&$\leq$&-&$(1,2,1,2,3,1)$\\
153&$>$&$>$&$\leq$&$>$&$>$&$\leq$&$\leq$&$>$&$>$&$>$&3(vi)&$(3,1,3,1,1,1)$\\
154&$>$&$>$&$\leq$&$>$&$>$&$\leq$&$\leq$&$>$&$>$&$\leq$&1&$(8^*,2)$\\
155&$>$&$>$&$\leq$&$>$&$>$&$\leq$&$\leq$&$>$&$\leq$&$>$&3(ix)&$(3,1,3,2,1)$ \\
156&$>$&$>$&$\leq$&$>$&$>$&$\leq$&$\leq$&$>$&$\leq$&$\leq$&1&$(7^*,3)$ \\
157&$>$&$>$&$\leq$&$>$&$>$&$\leq$&$\leq$&$\leq$&$>$&$>$&-&$(1,2,1,3,1,1,1)$ \\
158&$>$&$>$&$\leq$&$>$&$>$&$\leq$&$\leq$&$\leq$&$>$&$\leq$&1&$(8^*,2)$ \\
159&$>$&$>$&$\leq$&$>$&$>$&$\leq$&$\leq$&$\leq$&$\leq$&$>$&-&$(1,2,1,3,1,1,1)$ \\
160&$>$&$>$&$\leq$&$>$&$>$&$\leq$&$\leq$&$\leq$&$\leq$&$\leq$&-&$(1,2,1,3,1,1,1)$ \\
\end{tabular}
\newpage
 \begin{tabular}{ccccccccccccc}
 Case & A & B & C & D & E & F & G & H & I &J&Proposition& Inequalities \\ &&&&&&&&&&&&\\
161&$>$&$>$&$\leq$&$>$&$\leq$&$>$&$>$&$>$&$>$&$>$&3(ii)&$(1,2,1,2,1,1,1,1)$ \\
162&$>$&$>$&$\leq$&$>$&$\leq$&$>$&$>$&$>$&$>$&$\leq$&1&$(8^*,2)$ \\
163&$>$&$>$&$\leq$&$>$&$\leq$&$>$&$>$&$>$&$\leq$&$>$&3(iii)&$(1,2,2,1,1,2,1)$ \\
164&$>$&$>$&$\leq$&$>$&$\leq$&$>$&$>$&$>$&$\leq$&$\leq$&1&$(7^*,3)$ \\
165&$>$&$>$&$\leq$&$>$&$\leq$&$>$&$>$&$\leq$&$>$&$>$&3(iii)&$(1,2,2,1,2,1,1)$ \\
166&$>$&$>$&$\leq$&$>$&$\leq$&$>$&$>$&$\leq$&$>$&$\leq$&1&$(8^*,2)$ \\
167&$>$&$>$&$\leq$&$>$&$\leq$&$>$&$>$&$\leq$&$\leq$&$>$&3(ix)&$(3,2,1,3,1)$\\
168&$>$&$>$&$\leq$&$>$&$\leq$&$>$&$>$&$\leq$&$\leq$&$\leq$&-&$(1,2,2,1,3,1)$\\
169&$>$&$>$&$\leq$&$>$&$\leq$&$>$&$\leq$&$>$&$>$&$>$&3(iii)&$(1,2,2,2,1,1,1)$ \\
170&$>$&$>$&$\leq$&$>$&$\leq$&$>$&$\leq$&$>$&$>$&$\leq$&1&$(8^*,2)$ \\
171&$>$&$>$&$\leq$&$>$&$\leq$&$>$&$\leq$&$>$&$\leq$&$>$&3(iv)&$(1,2,2,2,2,1)$\\
172&$>$&$>$&$\leq$&$>$&$\leq$&$>$&$\leq$&$>$&$\leq$&$\leq$&1&$(7^*,3)$\\
173&$>$&$>$&$\leq$&$>$&$\leq$&$>$&$\leq$&$\leq$&$>$&$>$&3(ix)&$(3,2,3,1,1)$ \\
174&$>$&$>$&$\leq$&$>$&$\leq$&$>$&$\leq$&$\leq$&$>$&$\leq$&1&$(8^*,2)$ \\
175&$>$&$>$&$\leq$&$>$&$\leq$&$>$&$\leq$&$\leq$&$\leq$&$>$&-&$(1,2,2,3,1,1)$ \\
176&$>$&$>$&$\leq$&$>$&$\leq$&$>$&$\leq$&$\leq$&$\leq$&$\leq$&-&$(1,2,2,3,1,1)$ \\
177&$>$&$>$&$\leq$&$>$&$\leq$&$\leq$&$>$&$>$&$>$&$>$&3(vi)&$(3,3,1,1,1,1)$ \\
178&$>$&$>$&$\leq$&$>$&$\leq$&$\leq$&$>$&$>$&$>$&$\leq$&1&$(8^*,2)$ \\
179&$>$&$>$&$\leq$&$>$&$\leq$&$\leq$&$>$&$>$&$\leq$&$>$&3(ix)&$(3,3,1,2,1)$ \\
180&$>$&$>$&$\leq$&$>$&$\leq$&$\leq$&$>$&$>$&$\leq$&$\leq$&1&$(7^*,3)$ \\
181&$>$&$>$&$\leq$&$>$&$\leq$&$\leq$&$>$&$\leq$&$>$&$>$&3(ix)&$(3,3,2,1,1)$\\
182&$>$&$>$&$\leq$&$>$&$\leq$&$\leq$&$>$&$\leq$&$>$&$\leq$&1&$(8^*,2)$\\
183&$>$&$>$&$\leq$&$>$&$\leq$&$\leq$&$>$&$\leq$&$\leq$&$>$&3(vii)&$(3,3,3,1)$ \\
184&$>$&$>$&$\leq$&$>$&$\leq$&$\leq$&$>$&$\leq$&$\leq$&$\leq$&-&$(1,2,3,3,1)$ \\
185&$>$&$>$&$\leq$&$>$&$\leq$&$\leq$&$\leq$&$>$&$>$&$>$&-&$(1,2,3,1,1,1,1)$ \\
186&$>$&$>$&$\leq$&$>$&$\leq$&$\leq$&$\leq$&$>$&$>$&$\leq$&1&$(8^*,2)$ \\
187&$>$&$>$&$\leq$&$>$&$\leq$&$\leq$&$\leq$&$>$&$\leq$&$>$&-&$(1,2,3,1,2,1)$ \\
188&$>$&$>$&$\leq$&$>$&$\leq$&$\leq$&$\leq$&$>$&$\leq$&$\leq$&1&$(7^*,3)$ \\
189&$>$&$>$&$\leq$&$>$&$\leq$&$\leq$&$\leq$&$\leq$&$>$&$>$&6&$-$\\
190&$>$&$>$&$\leq$&$>$&$\leq$&$\leq$&$\leq$&$\leq$&$>$&$\leq$&1&$(8^*,2)$\\
191&$>$&$>$&$\leq$&$>$&$\leq$&$\leq$&$\leq$&$\leq$&$\leq$&$>$&-&$(3,3,1,1,1,1)$ \\
192&$>$&$>$&$\leq$&$>$&$\leq$&$\leq$&$\leq$&$\leq$&$\leq$&$\leq$&-&$(1,2,3,1,1,1,1)$ \\
193&$>$&$>$&$\leq$&$\leq$&$>$&$>$&$>$&$>$&$>$&$>$&3(v)&$(1,3,1,1,1,1,1,1)$ \\
194&$>$&$>$&$\leq$&$\leq$&$>$&$>$&$>$&$>$&$>$&$\leq$&1&$(8^*,2)$ \\
195&$>$&$>$&$\leq$&$\leq$&$>$&$>$&$>$&$>$&$\leq$&$>$&-&$(1,3,1,1,1,2,1)$ \\
196&$>$&$>$&$\leq$&$\leq$&$>$&$>$&$>$&$>$&$\leq$&$\leq$&1&$(7^*,3)$ \\
197&$>$&$>$&$\leq$&$\leq$&$>$&$>$&$>$&$\leq$&$>$&$>$&-&$(1,3,1,1,2,1,1)$\\
198&$>$&$>$&$\leq$&$\leq$&$>$&$>$&$>$&$\leq$&$>$&$\leq$&1&$(8^*,2)$\\
199&$>$&$>$&$\leq$&$\leq$&$>$&$>$&$>$&$\leq$&$\leq$&$>$&-&$(1,3,1,1,3,1)$ \\
200&$>$&$>$&$\leq$&$\leq$&$>$&$>$&$>$&$\leq$&$\leq$&$\leq$&-&$(1,2,1,1,1,2,1,1)$ \\
\end{tabular}
\newpage
 \begin{tabular}{ccccccccccccc}
 Case & A & B & C & D & E & F & G & H & I &J&Proposition& Inequalities \\ &&&&&&&&&&&&\\
201&$>$&$>$&$\leq$&$\leq$&$>$&$>$&$\leq$&$>$&$>$&$>$&-&$(1,3,1,2,1,1,1)$ \\
202&$>$&$>$&$\leq$&$\leq$&$>$&$>$&$\leq$&$>$&$>$&$\leq$&1&$(8^*,2)$ \\
203&$>$&$>$&$\leq$&$\leq$&$>$&$>$&$\leq$&$>$&$\leq$&$>$&-&$(1,3,1,2,2,1)$\\
204&$>$&$>$&$\leq$&$\leq$&$>$&$>$&$\leq$&$>$&$\leq$&$\leq$&1&$(7^*,3)$\\
205&$>$&$>$&$\leq$&$\leq$&$>$&$>$&$\leq$&$\leq$&$>$&$>$&-&$(1,3,1,3,1,1)$ \\
206&$>$&$>$&$\leq$&$\leq$&$>$&$>$&$\leq$&$\leq$&$>$&$\leq$&1&$(8^*,2)$ \\
207&$>$&$>$&$\leq$&$\leq$&$>$&$>$&$\leq$&$\leq$&$\leq$&$>$&-&$(1,2,1,1,2,1,1,1)$ \\
208&$>$&$>$&$\leq$&$\leq$&$>$&$>$&$\leq$&$\leq$&$\leq$&$\leq$&-&$(3,1,1,3,1,1)$ \\
209&$>$&$>$&$\leq$&$\leq$&$>$&$\leq$&$>$&$>$&$>$&$>$&-&$(1,3,2,1,1,1,1)$ \\
210&$>$&$>$&$\leq$&$\leq$&$>$&$\leq$&$>$&$>$&$>$&$\leq$&1&$(8^*,2)$ \\
211&$>$&$>$&$\leq$&$\leq$&$>$&$\leq$&$>$&$>$&$\leq$&$>$&-&$(1,3,2,1,2,1)$\\
212&$>$&$>$&$\leq$&$\leq$&$>$&$\leq$&$>$&$>$&$\leq$&$\leq$&1&$(7^*,3)$\\
213&$>$&$>$&$\leq$&$\leq$&$>$&$\leq$&$>$&$\leq$&$>$&$>$&-&$(1,3,2,2,1,1)$\\
214&$>$&$>$&$\leq$&$\leq$&$>$&$\leq$&$>$&$\leq$&$>$&$\leq$&1&$(8^*,2)$\\
215&$>$&$>$&$\leq$&$\leq$&$>$&$\leq$&$>$&$\leq$&$\leq$&$>$&-&$(1,3,2,3,1)$\\
216&$>$&$>$&$\leq$&$\leq$&$>$&$\leq$&$>$&$\leq$&$\leq$&$\leq$&4&$(3,1,\cdots,1),(1,2,1,2,2,1,1)$\\
217&$>$&$>$&$\leq$&$\leq$&$>$&$\leq$&$\leq$&$>$&$>$&$>$&-&$(1,3,3,1,1,1)$ \\
218&$>$&$>$&$\leq$&$\leq$&$>$&$\leq$&$\leq$&$>$&$>$&$\leq$&1&$(8^*,2)$ \\
219&$>$&$>$&$\leq$&$\leq$&$>$&$\leq$&$\leq$&$>$&$\leq$&$>$&-&$(1,3,3,2,1)$ \\
220&$>$&$>$&$\leq$&$\leq$&$>$&$\leq$&$\leq$&$>$&$\leq$&$\leq$&1&$(7^*,3)$ \\
221&$>$&$>$&$\leq$&$\leq$&$>$&$\leq$&$\leq$&$\leq$&$>$&$>$&-&$(1,3,3,1,1,1)$ \\
222&$>$&$>$&$\leq$&$\leq$&$>$&$\leq$&$\leq$&$\leq$&$>$&$\leq$&1&$(8^*,2)$ \\
223&$>$&$>$&$\leq$&$\leq$&$>$&$\leq$&$\leq$&$\leq$&$\leq$&$>$&-&$(1,2,1,2,2,1,1)$ \\
224&$>$&$>$&$\leq$&$\leq$&$>$&$\leq$&$\leq$&$\leq$&$\leq$&$\leq$&4&$(3,1,\cdots,1),(1,2,1,2,1,1,1,1)$ \\
225&$>$&$>$&$\leq$&$\leq$&$\leq$&$>$&$>$&$>$&$>$&$>$&-&$(1,2,2,1,1,1,1,1)$ \\
226&$>$&$>$&$\leq$&$\leq$&$\leq$&$>$&$>$&$>$&$>$&$\leq$&1&$(8^*,2)$ \\
227&$>$&$>$&$\leq$&$\leq$&$\leq$&$>$&$>$&$>$&$\leq$&$>$&-&$(1,3,1,1,2,1,1)$ \\
228&$>$&$>$&$\leq$&$\leq$&$\leq$&$>$&$>$&$>$&$\leq$&$\leq$&1&$(7^*,3)$ \\
229&$>$&$>$&$\leq$&$\leq$&$\leq$&$>$&$>$&$\leq$&$>$&$>$&-&$(1,3,1,1,2,1,1)$ \\
230&$>$&$>$&$\leq$&$\leq$&$\leq$&$>$&$>$&$\leq$&$>$&$\leq$&1&$(8^*,2)$ \\
231&$>$&$>$&$\leq$&$\leq$&$\leq$&$>$&$>$&$\leq$&$\leq$&$>$&-&$(1,2,1,1,1,2,1,1)$\\
232&$>$&$>$&$\leq$&$\leq$&$\leq$&$>$&$>$&$\leq$&$\leq$&$\leq$&-&$(1,2,1,1,1,4),(1,2,1,1,1,2,1,1)$\\
233&$>$&$>$&$\leq$&$\leq$&$\leq$&$>$&$\leq$&$>$&$>$&$>$&-&$(1,2,1,1,2,1,1,1)$ \\
234&$>$&$>$&$\leq$&$\leq$&$\leq$&$>$&$\leq$&$>$&$>$&$\leq$&1&$(8^*,2)$ \\
235&$>$&$>$&$\leq$&$\leq$&$\leq$&$>$&$\leq$&$>$&$\leq$&$>$&-&$(1,2,1,1,2,2,1)$ \\
236&$>$&$>$&$\leq$&$\leq$&$\leq$&$>$&$\leq$&$>$&$\leq$&$\leq$&1&$(7^*,3)$ \\
237&$>$&$>$&$\leq$&$\leq$&$\leq$&$>$&$\leq$&$\leq$&$>$&$>$&-&$(2,1,1,1,2,1,1,1)$ \\
238&$>$&$>$&$\leq$&$\leq$&$\leq$&$>$&$\leq$&$\leq$&$>$&$\leq$&1&$(8^*,2)$ \\
239&$>$&$>$&$\leq$&$\leq$&$\leq$&$>$&$\leq$&$\leq$&$\leq$&$>$&-&$(2,1,1,1,2,1,1,1)$ \\
240&$>$&$>$&$\leq$&$\leq$&$\leq$&$>$&$\leq$&$\leq$&$\leq$&$\leq$&4&$(3,1,\cdots,1)(1,2,1,1,2,1,1,1)$ \\
\end{tabular}
\newpage
 \begin{tabular}{ccccccccccccc}
 Case & A & B & C & D & E & F & G & H & I &J&Proposition& Inequalities \\ &&&&&&&&&&&&\\
241&$>$&$>$&$\leq$&$\leq$&$\leq$&$\leq$&$>$&$>$&$>$&$>$&-&$(1,2,2,1,1,1,1),(4,1,\cdots,1)$ \\
242&$>$&$>$&$\leq$&$\leq$&$\leq$&$\leq$&$>$&$>$&$>$&$\leq$&1&$(8^*,2)$ \\
243&$>$&$>$&$\leq$&$\leq$&$\leq$&$\leq$&$>$&$>$&$\leq$&$>$&-&$(1,2,1,\cdots,1)$ \\
244&$>$&$>$&$\leq$&$\leq$&$\leq$&$\leq$&$>$&$>$&$\leq$&$\leq$&1&$(7^*,3)$ \\
245&$>$&$>$&$\leq$&$\leq$&$\leq$&$\leq$&$>$&$\leq$&$>$&$>$&-&$(1,2,1,\cdots,1),(4,1,\cdots,1)$ \\
246&$>$&$>$&$\leq$&$\leq$&$\leq$&$\leq$&$>$&$\leq$&$>$&$\leq$&1&$(8^*,2)$ \\
247&$>$&$>$&$\leq$&$\leq$&$\leq$&$\leq$&$>$&$\leq$&$\leq$&$>$&-&$(3,1,1,1,2,1,1)$ \\
248&$>$&$>$&$\leq$&$\leq$&$\leq$&$\leq$&$>$&$\leq$&$\leq$&$\leq$&-&$(3,1,1,1,2,1,1),(3,1,\cdots,1)$ \\
249&$>$&$>$&$\leq$&$\leq$&$\leq$&$\leq$&$\leq$&$>$&$>$&$>$&6&$-$\\
250&$>$&$>$&$\leq$&$\leq$&$\leq$&$\leq$&$\leq$&$>$&$>$&$\leq$&1&$(8^*,2)$ \\
251&$>$&$>$&$\leq$&$\leq$&$\leq$&$\leq$&$\leq$&$>$&$\leq$&$>$&-&$(2,2,2,1,1,1,1)$ \\
252&$>$&$>$&$\leq$&$\leq$&$\leq$&$\leq$&$\leq$&$>$&$\leq$&$\leq$&1&$(7^*,3)$ \\
253&$>$&$>$&$\leq$&$\leq$&$\leq$&$\leq$&$\leq$&$\leq$&$>$&$>$&6&$-$ \\
254&$>$&$>$&$\leq$&$\leq$&$\leq$&$\leq$&$\leq$&$\leq$&$>$&$\leq$&1&$(8^*,2)$ \\
255&$>$&$>$&$\leq$&$\leq$&$\leq$&$\leq$&$\leq$&$\leq$&$\leq$&$>$&-&$(3,1,\cdots,1)$ \\
256&$>$&$>$&$\leq$&$\leq$&$\leq$&$\leq$&$\leq$&$\leq$&$\leq$&$\leq$&-&$(3,1,\cdots,1)$ \\
257&$>$&$\leq$&$>$&$>$&$>$&$>$&$>$&$>$&$>$&$>$&3(i)&$(2,1,1,1,1,1,1,1,1)$\\
258&$>$&$\leq$&$>$&$>$&$>$&$>$&$>$&$>$&$>$&$\leq$&1&$(8^*,2)$\\
259&$>$&$\leq$&$>$&$>$&$>$&$>$&$>$&$>$&$\leq$ &$>$&3(ii)&$(2,1,1,1,1,1,2,1)$ \\
260&$>$&$\leq$&$>$&$>$&$>$&$>$&$>$&$>$&$\leq$ &$\leq$&1&$(7^*,3)$ \\
261&$>$&$\leq$&$>$&$>$&$>$&$>$&$>$&$\leq$&$>$ &$>$&3(ii)&$(2,1,1,1,1,2,1,1)$\\
262&$>$&$\leq$&$>$&$>$&$>$&$>$&$>$&$\leq$&$>$ &$\leq$&1&$(8^*,2)$\\
263&$>$&$\leq$&$>$&$>$&$>$&$>$&$>$&$\leq$&$\leq$ &$>$&3(viii)&$(2,1,1,1,1,3,1)$ \\
264&$>$&$\leq$&$>$&$>$&$>$&$>$&$>$&$\leq$&$\leq$ &$\leq$&-&$(2,1,1,1,1,3,1)$ \\
265&$>$&$\leq$&$>$&$>$&$>$&$>$&$\leq$&$>$&$>$&$>$&3(ii)&$(2,1,1,1,2,1,1,1)$ \\
266&$>$&$\leq$&$>$&$>$&$>$&$>$&$\leq$&$>$&$>$&$\leq$&1&$(8^*,2)$ \\
267&$>$&$\leq$&$>$&$>$&$>$&$>$&$\leq$&$>$&$\leq$ &$>$&3(iii)&$(2,1,1,1,2,2,1)$ \\
268&$>$&$\leq$&$>$&$>$&$>$&$>$&$\leq$&$>$&$\leq$ &$\leq$&1&$(7^*,3)$ \\
269&$>$&$\leq$&$>$&$>$&$>$&$>$&$\leq$&$\leq$&$>$ &$>$&3(viii)&$(2,1,1,1,3,1,1)$ \\
270&$>$&$\leq$&$>$&$>$&$>$&$>$&$\leq$&$\leq$&$>$ &$\leq$&1&$(8^*,2)$ \\
271&$>$&$\leq$&$>$&$>$&$>$&$>$&$\leq$&$\leq$&$\leq$ &$>$&-&(2,1,1,1,3,1,1)\\
272&$>$&$\leq$&$>$&$>$&$>$&$>$&$\leq$&$\leq$&$\leq$ &$\leq$&5(ii)&$(2,1,\cdots,1),(2,1,1,1,3,1,1)$ \\
273&$>$&$\leq$&$>$&$>$&$>$&$\leq$&$>$&$>$&$>$ &$>$&3(ii)&$(2,1,1,2,1,1,1,1)$ \\
274&$>$&$\leq$&$>$&$>$&$>$&$\leq$&$>$&$>$&$>$ &$\leq$&1&$(8^*,2)$ \\
275&$>$&$\leq$&$>$&$>$&$>$&$\leq$&$>$&$>$&$\leq$ &$>$&3(iii)&$(2,1,1,2,1,2,1)$ \\
276&$>$&$\leq$&$>$&$>$&$>$&$\leq$&$>$&$>$&$\leq$ &$\leq$&1&$(7^*,3)$ \\
277&$>$&$\leq$&$>$&$>$&$>$&$\leq$&$>$&$\leq$&$>$ &$>$&3(iii)&(2,1,1,2,2,1,1) \\
278&$>$&$\leq$&$>$&$>$&$>$&$\leq$&$>$&$\leq$&$>$ &$\leq$&1& $(8^*,2)$\\
279&$>$&$\leq$&$>$&$>$&$>$&$\leq$&$>$&$\leq$&$\leq$ &$>$&3(ix)& $(2,1,1,2,3,1)$ \\
280&$>$&$\leq$&$>$&$>$&$>$&$\leq$&$>$&$\leq$&$\leq$ &$\leq$&5(ii)& $(2,1,\cdots,1),(2,1,1,2,3,1)$ \\
\end{tabular}
\newpage
 \begin{tabular}{ccccccccccccc}
 Case & A & B & C & D & E & F & G & H & I &J&Proposition& Inequalities \\ &&&&&&&&&&&&\\
281&$>$&$\leq$&$>$&$>$&$>$&$\leq$&$\leq$&$>$&$>$ &$>$&3(viii)& $(2,1,1,3,1,1,1)$\\
282&$>$&$\leq$&$>$&$>$&$>$&$\leq$&$\leq$&$>$&$>$ &$\leq$&1& $(8^*,2)$\\
283&$>$&$\leq$&$>$&$>$&$>$&$\leq$&$\leq$&$>$&$\leq$ &$>$&3(ix)& $(2,1,1,3,2,1)$ \\
284&$>$&$\leq$&$>$&$>$&$>$&$\leq$&$\leq$&$>$&$\leq$ &$\leq$&1& $(7^*,3)$ \\
285&$>$&$\leq$&$>$&$>$&$>$&$\leq$&$\leq$&$\leq$&$>$ &$>$&-&(2,1,1,3,1,1,1)\\
286&$>$&$\leq$&$>$&$>$&$>$&$\leq$&$\leq$&$\leq$&$>$ &$\leq$&1&$(8^*,2)$\\
287&$>$&$\leq$&$>$&$>$&$>$&$\leq$&$\leq$&$\leq$&$\leq$ &$>$&5(ii)& $(2,1,\cdots,1),(2,1,1,3,1,1,1)$\\
288&$>$&$\leq$&$>$&$>$&$>$&$\leq$&$\leq$&$\leq$&$\leq$ &$\leq$&5(i)& $(2,1,\cdots,1),(2,1,1,3,1,1,1)$\\
289&$>$&$\leq$&$>$&$>$&$\leq$&$>$&$>$&$>$&$>$ &$>$&3(ii)&$(2,1,2,1,1,1,1,1)$ \\
290&$>$&$\leq$&$>$&$>$&$\leq$&$>$&$>$&$>$&$>$ &$\leq$&1&$(8^*,2)$ \\
291&$>$&$\leq$&$>$&$>$&$\leq$&$>$&$>$&$>$&$\leq$ &$>$&3(iii)& $(2,1,2,1,1,2,1)$ \\
292&$>$&$\leq$&$>$&$>$&$\leq$&$>$&$>$&$>$&$\leq$ &$\leq$&1& $(7^*,3)$ \\
293&$>$&$\leq$&$>$&$>$&$\leq$&$>$&$>$&$\leq$&$>$ &$>$&3(iii)&$(2,1,2,1,2,1,1)$\\
294&$>$&$\leq$&$>$&$>$&$\leq$&$>$&$>$&$\leq$&$>$ &$\leq$&1&$(8^*,2)$\\
295&$>$&$\leq$&$>$&$>$&$\leq$&$>$&$>$&$\leq$&$\leq$ &$>$&3(ix)&$(2,1,2,1,3,1)$ \\
296&$>$&$\leq$&$>$&$>$&$\leq$&$>$&$>$&$\leq$&$\leq$ &$\leq$&5(ii)&$(2,1,\cdots,1),(2,1,2,1,3,1)$ \\
297&$>$&$\leq$&$>$&$>$&$\leq$&$>$&$\leq$&$>$&$>$ &$>$&3(iii)&$(2,1,2,2,1,1,1)$\\
298&$>$&$\leq$&$>$&$>$&$\leq$&$>$&$\leq$&$>$&$>$ &$\leq$&1&$(8^*,2)$\\
299&$>$&$\leq$&$>$&$>$&$\leq$&$>$&$\leq$&$>$&$\leq$ &$>$&3(iv)& $(2,1,2,2,2,1)$ \\
300&$>$&$\leq$&$>$&$>$&$\leq$&$>$&$\leq$&$>$&$\leq$ &$\leq$&1& $(7^*,3)$ \\
301&$>$&$\leq$&$>$&$>$&$\leq$&$>$&$\leq$&$\leq$&$>$ &$>$&3(ix)&$(2,1,2,3,1,1)$ \\
302&$>$&$\leq$&$>$&$>$&$\leq$&$>$&$\leq$&$\leq$&$>$ &$\leq$&1&$(8^*,2)$ \\
303&$>$&$\leq$&$>$&$>$&$\leq$&$>$&$\leq$&$\leq$&$\leq$ &$>$&5(ii)& $(2,1,\cdots,1),(2,1,2,3,1,1)$\\
304&$>$&$\leq$&$>$&$>$&$\leq$&$>$&$\leq$&$\leq$&$\leq$ &$\leq$&5(i)&$(2,1,\cdots,1),(2,1,2,3,1,1)$\\
305&$>$&$\leq$&$>$&$>$&$\leq$&$\leq$&$>$&$>$&$>$ &$>$&3(viii)&$(2,1,3,1,1,1,1)$ \\
306&$>$&$\leq$&$>$&$>$&$\leq$&$\leq$&$>$&$>$&$>$ &$\leq$&1&$(8^*,2)$ \\
307&$>$&$\leq$&$>$&$>$&$\leq$&$\leq$&$>$&$>$&$\leq$ &$>$&3(ix)&$(2,1,3,1,2,1)$\\
308&$>$&$\leq$&$>$&$>$&$\leq$&$\leq$&$>$&$>$&$\leq$ &$\leq$&1&$(7^*,3)$\\
309&$>$&$\leq$&$>$&$>$&$\leq$&$\leq$&$>$&$\leq$&$>$ &$>$&3(ix)&$(2,1,3,2,1,1)$\\
310&$>$&$\leq$&$>$&$>$&$\leq$&$\leq$&$>$&$\leq$&$>$ &$\leq$&1&$(8^*,2)$\\
311&$>$&$\leq$&$>$&$>$&$\leq$&$\leq$&$>$&$\leq$&$\leq$&$>$&3(xi)&$(2,1,3,3,1)$\\
312&$>$&$\leq$&$>$&$>$&$\leq$&$\leq$&$>$&$\leq$&$\leq$&$\leq$&5(i)&$(2,1,\cdots,1),(2,1,3,3,1)$\\
313&$>$&$\leq$&$>$&$>$&$\leq$&$\leq$&$\leq$&$>$&$>$ &$>$&-&$(2,1,3,1,1,1,1)$\\
314&$>$&$\leq$&$>$&$>$&$\leq$&$\leq$&$\leq$&$>$&$>$ &$\leq$&1&$(8^*,2)$\\
315&$>$&$\leq$&$>$&$>$&$\leq$&$\leq$&$\leq$&$>$&$\leq$ &$>$&5(ii)&$(2,1,\cdots,1),(2,1,3,1,2,1)$ \\
316&$>$&$\leq$&$>$&$>$&$\leq$&$\leq$&$\leq$&$>$&$\leq$ &$\leq$&1&$(7^*,3)$ \\
317&$>$&$\leq$&$>$&$>$&$\leq$&$\leq$&$\leq$&$\leq$&$>$ & $>$&5(ii)&$(2,1,\cdots,1),(2,1,3,1,1,1,1)$\\
318&$>$&$\leq$&$>$&$>$&$\leq$&$\leq$&$\leq$&$\leq$&$>$ & $\leq$&1&$(8^*,2)$\\
319&$>$&$\leq$&$>$&$>$&$\leq$&$\leq$&$\leq$&$\leq$&$\leq$&$>$&5(i)&$(2,1,\cdots,1),(2,1,3,1,1,1,1)$\\
320&$>$&$\leq$&$>$&$>$&$\leq$&$\leq$&$\leq$&$\leq$&$\leq$&$\leq$&2&$(2,1,1,1,1,1,1,1,1)$\\
\end{tabular}
\newpage
 \begin{tabular}{ccccccccccccc}
 Case & A & B & C & D & E & F & G & H & I &J&Proposition& Inequalities \\ &&&&&&&&&&&&\\
321&$>$&$\leq$&$>$&$\leq$&$>$&$>$&$>$&$>$&$>$&$>$&3(ii)&$(2,2,1,1,1,1,1,1)$\\
322&$>$&$\leq$&$>$&$\leq$&$>$&$>$&$>$&$>$&$>$&$\leq$&1&$(8^*,2)$\\
323&$>$&$\leq$&$>$&$\leq$&$>$&$>$&$>$&$>$&$\leq$&$>$&3(iii)&$(2,2,1,1,1,2,1)$ \\
324&$>$&$\leq$&$>$&$\leq$&$>$&$>$&$>$&$>$&$\leq$&$\leq$&1&$(7^*,3)$ \\
325&$>$&$\leq$&$>$&$\leq$&$>$&$>$&$>$&$\leq$&$>$&$>$&3(iii)&$(2,2,1,1,2,1,1)$\\
326&$>$&$\leq$&$>$&$\leq$&$>$&$>$&$>$&$\leq$&$>$&$\leq$&1&$(8^*,2)$\\
327&$>$&$\leq$&$>$&$\leq$&$>$&$>$&$>$&$\leq$&$\leq$&$>$&3(ix)&$(2,2,1,1,3,1)$\\
328&$>$&$\leq$&$>$&$\leq$&$>$&$>$&$>$&$\leq$&$\leq$&$\leq$&5(ii)&$(2,1,\cdots,1),(2,2,1,1,3,1)$\\
329&$>$&$\leq$&$>$&$\leq$&$>$&$>$&$\leq$&$>$&$>$&$>$&3(iii)&$(2,2,1,2,1,1,1)$\\
330&$>$&$\leq$&$>$&$\leq$&$>$&$>$&$\leq$&$>$&$>$&$\leq$&1&$(8^*,2)$\\
331&$>$&$\leq$&$>$&$\leq$&$>$&$>$&$\leq$&$>$&$\leq$&$>$&3(iv)&$(2,2,1,2,2,1)$\\
332&$>$&$\leq$&$>$&$\leq$&$>$&$>$&$\leq$&$>$&$\leq$&$\leq$&1&$(7^*,3)$\\
333&$>$&$\leq$&$>$&$\leq$&$>$&$>$&$\leq$&$\leq$&$>$&$>$&-&$(2,2,1,3,1,1)$\\
334&$>$&$\leq$&$>$&$\leq$&$>$&$>$&$\leq$&$\leq$&$>$&$\leq$&1&$(8^*,2)$\\
335&$>$&$\leq$&$>$&$\leq$&$>$&$>$&$\leq$&$\leq$&$\leq$&$>$&5(ii)&$(2,1,\cdots,1),(2,2,1,3,1,1)$\\
336&$>$&$\leq$&$>$&$\leq$&$>$&$>$&$\leq$&$\leq$&$\leq$&$\leq$&5(i)&$(2,1,\cdots,1),(2,2,1,3,1,1)$\\
337&$>$&$\leq$&$>$&$\leq$&$>$&$\leq$&$>$&$>$&$>$&$>$&3(iii)&$(2,2,2,1,1,1,1)$\\
338&$>$&$\leq$&$>$&$\leq$&$>$&$\leq$&$>$&$>$&$>$&$\leq$&1&$(8^*,2)$\\
339&$>$&$\leq$&$>$&$\leq$&$>$&$\leq$&$>$&$>$&$\leq$&$>$&3(iv)&$(2,2,2,1,2,1)$\\
340&$>$&$\leq$&$>$&$\leq$&$>$&$\leq$&$>$&$>$&$\leq$&$\leq$&1&$(7^*,3)$\\
341&$>$ &$\leq$&$>$&$\leq$&$>$&$\leq$&$>$&$\leq$&$>$&$>$&3(iv)&$(2,2,2,2,1,1)$\\
342&$>$ &$\leq$&$>$&$\leq$&$>$&$\leq$&$>$&$\leq$&$>$&$\leq$&1&$(8^*,2)$\\
343&$>$&$\leq$&$>$&$\leq$&$>$&$\leq$&$>$&$\leq$&$\leq$&$>$&3(x)&$(2,2,2,3,1)$\\
344&$>$&$\leq$&$>$&$\leq$&$>$&$\leq$&$>$&$\leq$&$\leq$&$\leq$&5(i)&$(2,1,\cdots,1),(2,2,2,3,1)$\\
345&$>$&$\leq$&$>$&$\leq$&$>$&$\leq$&$\leq$&$>$&$>$&$>$&3(ix)&$(2,2,3,1,1,1)$\\
346&$>$&$\leq$&$>$&$\leq$&$>$&$\leq$&$\leq$&$>$&$>$&$\leq$&1&$(8^*,2)$\\
347&$>$&$\leq$&$>$&$\leq$&$>$&$\leq$&$\leq$&$>$&$\leq$&$>$&3(x)&$(2,2,3,2,1)$\\
348&$>$&$\leq$&$>$&$\leq$&$>$&$\leq$&$\leq$&$>$&$\leq$&$\leq$&1&$(7^*,3)$\\
349&$>$&$\leq$&$>$&$\leq$&$>$&$\leq$&$\leq$&$\leq$&$>$&$>$&5(ii)& $(2,1,\cdots,1),(2,2,2,1,1,1,1)$\\
350&$>$&$\leq$&$>$&$\leq$&$>$&$\leq$&$\leq$&$\leq$&$>$&$\leq$&1& $(8^*,2)$\\
351&$>$&$\leq$&$>$&$\leq$&$>$&$\leq$&$\leq$&$\leq$&$\leq$&$>$&5(i)&$(2,1,\cdots,1),(2,2,2,1,1,1,1)$\\
352&$>$&$\leq$&$>$&$\leq$&$>$&$\leq$&$\leq$&$\leq$&$\leq$&$\leq$&2&$(2,1,1,1,1,1,1,1,1)$\\
353&$>$&$\leq$&$>$&$\leq$&$\leq$&$>$&$>$&$>$&$>$&$>$&3(viii)&$(2,3,1,1,1,1,1)$\\
354&$>$&$\leq$&$>$&$\leq$&$\leq$&$>$&$>$&$>$&$>$&$\leq$&1&$(8^*,2)$\\
355&$>$&$\leq$&$>$&$\leq$&$\leq$&$>$&$>$&$>$&$\leq$&$>$&3(ix)&$(2,3,1,1,2,1)$\\
356&$>$&$\leq$&$>$&$\leq$&$\leq$&$>$&$>$&$>$&$\leq$&$\leq$&1&$(7^*,3)$\\
357&$>$&$\leq$&$>$&$\leq$&$\leq$&$>$&$>$&$\leq$&$>$&$>$&3(ix)&$(2,3,1,2,1,1)$\\
358&$>$&$\leq$&$>$&$\leq$&$\leq$&$>$&$>$&$\leq$&$>$&$\leq$&1&$(8^*,2)$\\
359&$>$ &$\leq$&$>$&$\leq$&$\leq$&$>$&$>$&$\leq$&$\leq$&$>$&3(xi)&$(2,3,1,3,1)$ \\
360&$>$ &$\leq$&$>$&$\leq$&$\leq$&$>$&$>$&$\leq$&$\leq$&$\leq$&5(i)&$(2,1,\cdots,1),(2,3,1,3,1)$ \\
\end{tabular}
\newpage
 \begin{tabular}{ccccccccccccc}
 Case & A & B & C & D & E & F & G & H & I &J&Proposition& Inequalities \\ &&&&&&&&&&&&\\
361&$>$&$\leq$&$>$&$\leq$&$\leq$&$>$&$\leq$&$>$&$>$&$>$&3(ix)& $(2,3,2,1,1,1)$\\
362&$>$&$\leq$&$>$&$\leq$&$\leq$&$>$&$\leq$&$>$&$>$&$\leq$&1& $(8^*,2)$\\
363&$>$&$\leq$&$>$&$\leq$&$\leq$&$>$&$\leq$&$>$&$\leq$&$>$&3(x)&$(2,3,2,2,1)$\\
364&$>$&$\leq$&$>$&$\leq$&$\leq$&$>$&$\leq$&$>$&$\leq$&$\leq$&1&$(7^*,3)$\\
365&$>$&$\leq$&$>$&$\leq$&$\leq$&$>$&$\leq$&$\leq$&$>$&$>$&3(xi)&$(2,3,3,1,1)$\\
366&$>$&$\leq$&$>$&$\leq$&$\leq$&$>$&$\leq$&$\leq$&$>$&$\leq$&1&$(8^*,2)$\\
367&$>$&$\leq$&$>$&$\leq$&$\leq$&$>$&$\leq$&$\leq$&$\leq$&$>$&5(i)&$(2,1,\cdots,1),(2,3,3,1,1)$\\
368&$>$&$\leq$&$>$&$\leq$&$\leq$&$>$&$\leq$&$\leq$&$\leq$&$\leq$&2&$(2,1,1,1,1,1,1,1,1)$\\
369&$>$&$\leq$&$>$&$\leq$&$\leq$&$\leq$&$>$&$>$&$>$&$>$&-&$(2,3,1,1,1,1,1)$\\
370&$>$&$\leq$&$>$&$\leq$&$\leq$&$\leq$&$>$&$>$&$>$&$\leq$&1&$(8^*,2)$\\
371&$>$&$\leq$&$>$&$\leq$&$\leq$&$\leq$&$>$&$>$&$\leq$&$>$&5(ii)&$(2,1,\cdots,1),(2,3,1,1,2,1)$\\
372&$>$&$\leq$&$>$&$\leq$&$\leq$&$\leq$&$>$&$>$&$\leq$&$\leq$&1&$(7^*,3)$\\
373&$>$ &$\leq$&$>$&$\leq$&$\leq$&$\leq$&$>$&$\leq$&$>$&$>$&5(ii)&$(2,1,\cdots,1),(2,3,1,2,1,1)$\\
374&$>$ &$\leq$&$>$&$\leq$&$\leq$&$\leq$&$>$&$\leq$&$>$&$\leq$&1&$(8^*,2)$\\
375&$>$&$\leq$&$>$&$\leq$&$\leq$&$\leq$&$>$&$\leq$&$\leq$&$>$&5(i)&$(2,1,\cdots,1),(2,3,1,1,1,1,1)$\\
376&$>$&$\leq$&$>$&$\leq$&$\leq$&$\leq$&$>$&$\leq$&$\leq$&$\leq$&2&$(2,1,1,1,1,1,1,1,1)$\\
377&$>$&$\leq$&$>$&$\leq$&$\leq$&$\leq$&$\leq$&$>$&$>$&$>$&5(ii)&$(2,1,\cdots,1),(2,3,1,1,1,1,1)$\\
378&$>$&$\leq$&$>$&$\leq$&$\leq$&$\leq$&$\leq$&$>$&$>$&$\leq$&1&$(8^*,2)$,\\
379&$>$&$\leq$&$>$&$\leq$&$\leq$&$\leq$&$\leq$&$>$&$\leq$&$>$&5(i)&$(2,1,\cdots,1),(2,3,1,1,1,1,1)$\\
380&$>$&$\leq$&$>$&$\leq$&$\leq$&$\leq$&$\leq$&$>$&$\leq$&$\leq$&1&$(7^*,3)$\\
381&$>$&$\leq$&$>$&$\leq$&$\leq$&$\leq$&$\leq$&$\leq$&$>$ &$>$&5(i)&$(2,1,\cdots,1),(2,3,1,1,1,1,1)$\\
382&$>$&$\leq$&$>$&$\leq$&$\leq$&$\leq$&$\leq$&$\leq$&$>$ &$\leq$&1&$(8^*,2)$\\
383&$>$&$\leq$&$>$&$\leq$&$\leq$&$\leq$&$\leq$&$\leq$&$\leq$&$>$&-&$(2,1,1,1,1,1,1,1,1)$\\
384&$>$&$\leq$&$>$&$\leq$&$\leq$&$\leq$&$\leq$&$\leq$&$\leq$&$\leq$&2&$(2,1,1,1,1,1,1,1,1)$\\
385&$>$&$\leq$&$\leq$&$>$&$>$&$>$&$>$&$>$&$>$&$>$&3(v)&$(3,1,1,1,1,1,1,1)$ \\
386&$>$&$\leq$&$\leq$&$>$&$>$&$>$&$>$&$>$&$>$&$\leq$&1&$(8^*,2)$ \\
387&$>$&$\leq$&$\leq$&$>$&$>$&$>$&$>$&$>$&$\leq$&$>$&3(viii)&$(3,1,1,1,1,2,1)$ \\
388&$>$&$\leq$&$\leq$&$>$&$>$&$>$&$>$&$>$&$\leq$&$\leq$&1&$(7^*,3)$ \\
389&$>$&$\leq$&$\leq$&$>$&$>$&$>$&$>$&$\leq$&$>$&$>$&3(viii)&$(3,1,1,1,2,1,1)$ \\
390&$>$&$\leq$&$\leq$&$>$&$>$&$>$&$>$&$\leq$&$>$&$\leq$&1&$(8^*,2)$ \\
391&$>$&$\leq$&$\leq$&$>$&$>$&$>$&$>$&$\leq$&$\leq$&$>$&3(vi)&$(3,1,1,1,3,1)$ \\
392&$>$&$\leq$&$\leq$&$>$&$>$&$>$&$>$&$\leq$&$\leq$&$\leq$&5(ii)&$(2,1,\cdots,1),(3,1,1,1,3,1)$ \\
393&$>$&$\leq$&$\leq$&$>$&$>$&$>$&$\leq$&$>$&$>$&$>$&3(viii)&$(3,1,1,2,1,1,1)$ \\
394&$>$&$\leq$&$\leq$&$>$&$>$&$>$&$\leq$&$>$&$>$&$\leq$&1&$(8^*,2)$ \\
395&$>$&$\leq$&$\leq$&$>$&$>$&$>$&$\leq$&$>$&$\leq$&$>$&3(ix)&$(3,1,1,2,2,1)$ \\
396&$>$&$\leq$&$\leq$&$>$&$>$&$>$&$\leq$&$>$&$\leq$&$\leq$&1&$(7^*,3)$ \\
397&$>$&$\leq$&$\leq$&$>$&$>$&$>$&$\leq$&$\leq$&$>$&$>$&3(vi)&$(3,1,1,3,1,1)$ \\
398&$>$&$\leq$&$\leq$&$>$&$>$&$>$&$\leq$&$\leq$&$>$&$\leq$&1&$(8^*,2)$ \\
399&$>$&$\leq$&$\leq$&$>$&$>$&$>$&$\leq$&$\leq$&$\leq$&$>$&5(ii)&$(2,1,\cdots,1),(3,1,1,3,1,1)$ \\
400&$>$&$\leq$&$\leq$&$>$&$>$&$>$&$\leq$&$\leq$&$\leq$&$\leq$&5(i)&$(2,1,\cdots,1),(3,1,1,3,1,1)$ \\
\end{tabular}
\newpage
 \begin{tabular}{ccccccccccccc}
 Case & A & B & C & D & E & F & G & H & I &J&Proposition& Inequalities \\ &&&&&&&&&&&&\\
401&$>$&$\leq$&$\leq$&$>$&$>$&$\leq$&$>$&$>$&$>$&$>$&3(viii)&$(3,1,2,1,1,1,1)$ \\
402&$>$&$\leq$&$\leq$&$>$&$>$&$\leq$&$>$&$>$&$>$&$\leq$&1&$(8^*,2)$ \\
403&$>$&$\leq$&$\leq$&$>$&$>$&$\leq$&$>$&$>$&$\leq$&$>$&3(ix)&$(3,1,2,1,2,1)$ \\
404&$>$&$\leq$&$\leq$&$>$&$>$&$\leq$&$>$&$>$&$\leq$&$\leq$&1&$(7^*,3)$ \\
405&$>$&$\leq$&$\leq$&$>$&$>$&$\leq$&$>$&$\leq$&$>$&$>$&3(ix)&$(3,1,2,2,1,1)$ \\
406&$>$&$\leq$&$\leq$&$>$&$>$&$\leq$&$>$&$\leq$&$>$&$\leq$&1&$(8^*,2)$ \\
407&$>$&$\leq$&$\leq$&$>$&$>$&$\leq$&$>$&$\leq$&$\leq$&$>$&3(xi)&$(3,1,2,3,1)$\\
408&$>$&$\leq$&$\leq$&$>$&$>$&$\leq$&$>$&$\leq$&$\leq$&$\leq$&5(i)&$(2,1,\cdots,1),(3,1,2,3,1)$\\
409&$>$&$\leq$&$\leq$&$>$&$>$&$\leq$&$\leq$&$>$&$>$&$>$&3(vi)&$(3,1,3,1,1,1)$\\
410&$>$&$\leq$&$\leq$&$>$&$>$&$\leq$&$\leq$&$>$&$>$&$\leq$&1&$(8^*,2)$\\
411&$>$&$\leq$&$\leq$&$>$&$>$&$\leq$&$\leq$&$>$&$\leq$&$>$&3(xi)&$(3,1,3,2,1)$ \\
412&$>$&$\leq$&$\leq$&$>$&$>$&$\leq$&$\leq$&$>$&$\leq$&$\leq$&1&$(7^*,3)$ \\
413&$>$&$\leq$&$\leq$&$>$&$>$&$\leq$&$\leq$&$\leq$&$>$&$>$&5(ii)&$(2,1,\cdots,1),(3,1,3,1,1,1)$ \\
414&$>$&$\leq$&$\leq$&$>$&$>$&$\leq$&$\leq$&$\leq$&$>$&$\leq$&1&$(8^*,2)$ \\
415&$>$&$\leq$&$\leq$&$>$&$>$&$\leq$&$\leq$&$\leq$&$\leq$&$>$&5(i)&$(2,1,\cdots,1),(2,1,1,3,1,1,1)$ \\
416&$>$&$\leq$&$\leq$&$>$&$>$&$\leq$&$\leq$&$\leq$&$\leq$&$\leq$&2&(2,1,1,1,1,1,1,1,1) \\
417&$>$&$\leq$&$\leq$&$>$&$\leq$&$>$&$>$&$>$&$>$&$>$&3(viii)&$(3,2,1,1,1,1,1)$ \\
418&$>$&$\leq$&$\leq$&$>$&$\leq$&$>$&$>$&$>$&$>$&$\leq$&1&$(8^*,2)$ \\
419&$>$&$\leq$&$\leq$&$>$&$\leq$&$>$&$>$&$>$&$\leq$&$>$&3(ix)&$(3,2,1,1,2,1)$ \\
420&$>$&$\leq$&$\leq$&$>$&$\leq$&$>$&$>$&$>$&$\leq$&$\leq$&1&$(7^*,3)$ \\
421&$>$&$\leq$&$\leq$&$>$&$\leq$&$>$&$>$&$\leq$&$>$&$>$&3(ix)&$(3,2,1,2,1,1)$ \\
422&$>$&$\leq$&$\leq$&$>$&$\leq$&$>$&$>$&$\leq$&$>$&$\leq$&1&$(8^*,2)$ \\
423&$>$&$\leq$&$\leq$&$>$&$\leq$&$>$&$>$&$\leq$&$\leq$&$>$&3(xi)&$(3,2,1,3,1)$\\
424&$>$&$\leq$&$\leq$&$>$&$\leq$&$>$&$>$&$\leq$&$\leq$&$\leq$&-&$(3,2,1,3,1)$\\
425&$>$&$\leq$&$\leq$&$>$&$\leq$&$>$&$\leq$&$>$&$>$&$>$&3(ix)&$(3,2,2,1,1,1)$ \\
426&$>$&$\leq$&$\leq$&$>$&$\leq$&$>$&$\leq$&$>$&$>$&$\leq$&1&$(8^*,2)$ \\
427&$>$&$\leq$&$\leq$&$>$&$\leq$&$>$&$\leq$&$>$&$\leq$&$>$&3(x)&$(3,2,2,2,1)$\\
428&$>$&$\leq$&$\leq$&$>$&$\leq$&$>$&$\leq$&$>$&$\leq$&$\leq$&1&$(7^*,3)$\\
429&$>$&$\leq$&$\leq$&$>$&$\leq$&$>$&$\leq$&$\leq$&$>$&$>$&3(xi)&(3,2,3,1,1) \\
430&$>$&$\leq$&$\leq$&$>$&$\leq$&$>$&$\leq$&$\leq$&$>$&$\leq$&1&$(8^*,2)$ \\
431&$>$&$\leq$&$\leq$&$>$&$\leq$&$>$&$\leq$&$\leq$&$\leq$&$>$&5(i)&$(2,1,\cdots,1),(3,2,3,1,1)$ \\
432&$>$&$\leq$&$\leq$&$>$&$\leq$&$>$&$\leq$&$\leq$&$\leq$&$\leq$&2&$(2,1,1,1,1,1,1,1,1)$ \\
433&$>$&$\leq$&$\leq$&$>$&$\leq$&$\leq$&$>$&$>$&$>$&$>$&3(vi)&$(3,3,1,1,1,1)$ \\
434&$>$&$\leq$&$\leq$&$>$&$\leq$&$\leq$&$>$&$>$&$>$&$\leq$&1&$(8^*,2)$ \\
435&$>$&$\leq$&$\leq$&$>$&$\leq$&$\leq$&$>$&$>$&$\leq$&$>$&3(xi)&$(3,3,1,2,1)$ \\
436&$>$&$\leq$&$\leq$&$>$&$\leq$&$\leq$&$>$&$>$&$\leq$&$\leq$&1&$(7^*,3)$ \\
437&$>$&$\leq$&$\leq$&$>$&$\leq$&$\leq$&$>$&$\leq$&$>$&$>$&3(xi)&$(3,3,2,1,1)$\\
438&$>$&$\leq$&$\leq$&$>$&$\leq$&$\leq$&$>$&$\leq$&$>$&$\leq$&1&$(8^*,2)$\\
439&$>$&$\leq$&$\leq$&$>$&$\leq$&$\leq$&$>$&$\leq$&$\leq$&$>$&3(vii)&$(3,3,3,1)$ \\
440&$>$&$\leq$&$\leq$&$>$&$\leq$&$\leq$&$>$&$\leq$&$\leq$&$\leq$&2&$(2,1,1,1,1,1,1,1,1)$ \\
\end{tabular}
\newpage
 \begin{tabular}{ccccccccccccc}
 Case & A & B & C & D & E & F & G & H & I &J&Proposition& Inequalities \\ &&&&&&&&&&&&\\
441&$>$&$\leq$&$\leq$&$>$&$\leq$&$\leq$&$\leq$&$>$&$>$&$>$&5(ii)&$(2,1,\cdots,1),(3,3,1,1,1,1)$ \\
442&$>$&$\leq$&$\leq$&$>$&$\leq$&$\leq$&$\leq$&$>$&$>$&$\leq$&1&$(8^*,2)$ \\
443&$>$&$\leq$&$\leq$&$>$&$\leq$&$\leq$&$\leq$&$>$&$\leq$&$>$&-&$(3,2,1,1,2,1)$ \\
444&$>$&$\leq$&$\leq$&$>$&$\leq$&$\leq$&$\leq$&$>$&$\leq$&$\leq$&1&$(7^*,3)$ \\
445&$>$&$\leq$&$\leq$&$>$&$\leq$&$\leq$&$\leq$&$\leq$&$>$&$>$&-&$(2,1,2,1,1,1,1,1),(4,1,\cdots,1)$\\
446&$>$&$\leq$&$\leq$&$>$&$\leq$&$\leq$&$\leq$&$\leq$&$>$&$\leq$&1&$(8^*,2)$\\
447&$>$&$\leq$&$\leq$&$>$&$\leq$&$\leq$&$\leq$&$\leq$&$\leq$&$>$&2&$(2,1,1,1,1,1,1,1,1)$ \\
448&$>$&$\leq$&$\leq$&$>$&$\leq$&$\leq$&$\leq$&$\leq$&$\leq$&$\leq$&2&$(2,1,1,1,1,1,1,1,1)$ \\
449&$>$&$\leq$&$\leq$&$\leq$&$>$&$>$&$>$&$>$&$>$&$>$&-&$(2,1,\cdots,1),(4,1,\cdots,1)$ \\
450&$>$&$\leq$&$\leq$&$\leq$&$>$&$>$&$>$&$>$&$>$&$\leq$&1&$(8^*,2)$ \\
451&$>$&$\leq$&$\leq$&$\leq$&$>$&$>$&$>$&$>$&$\leq$&$>$&-&$(3,1,1,1,1,2,1)$ \\
452&$>$&$\leq$&$\leq$&$\leq$&$>$&$>$&$>$&$>$&$\leq$&$\leq$&1&$(7^*,3)$ \\
453&$>$&$\leq$&$\leq$&$\leq$&$>$&$>$&$>$&$\leq$&$>$&$>$&-&$(3,1,1,1,2,1,1)$\\
454&$>$&$\leq$&$\leq$&$\leq$&$>$&$>$&$>$&$\leq$&$>$&$\leq$&1&$(8^*,2)$\\
455&$>$&$\leq$&$\leq$&$\leq$&$>$&$>$&$>$&$\leq$&$\leq$&$>$&5(ii)&$(2,1,\cdots,1),(3,1,1,1,2,1,1)$ \\
456&$>$&$\leq$&$\leq$&$\leq$&$>$&$>$&$>$&$\leq$&$\leq$&$\leq$&5(i)&$(2,1,\cdots,1),(4,1,\cdots,1)$ \\
457&$>$&$\leq$&$\leq$&$\leq$&$>$&$>$&$\leq$&$>$&$>$&$>$&-&$(3,1,1,2,1,1,1)$ \\
458&$>$&$\leq$&$\leq$&$\leq$&$>$&$>$&$\leq$&$>$&$>$&$\leq$&1&$(8^*,2)$ \\
459&$>$&$\leq$&$\leq$&$\leq$&$>$&$>$&$\leq$&$>$&$\leq$&$>$&5(ii)&$(2,1,\cdots,1),(3,1,1,2,1,1,1)$\\
460&$>$&$\leq$&$\leq$&$\leq$&$>$&$>$&$\leq$&$>$&$\leq$&$\leq$&1&$(7^*,3)$\\
461&$>$&$\leq$&$\leq$&$\leq$&$>$&$>$&$\leq$&$\leq$&$>$&$>$&5(ii)&$(2,1,\cdots,1),(3,1,1,2,1,1,1)$ \\
462&$>$&$\leq$&$\leq$&$\leq$&$>$&$>$&$\leq$&$\leq$&$>$&$\leq$&1&$(8^*,2)$ \\
463&$>$&$\leq$&$\leq$&$\leq$&$>$&$>$&$\leq$&$\leq$&$\leq$&$>$&5(i)&$(2,1,\cdots,1),(4,1,\cdots,1)$ \\
464&$>$&$\leq$&$\leq$&$\leq$&$>$&$>$&$\leq$&$\leq$&$\leq$&$\leq$&2&$(2,1,1,1,1,1,1,1,1)$ \\
465&$>$&$\leq$&$\leq$&$\leq$&$>$&$\leq$&$>$&$>$&$>$&$>$&-&$(3,1,2,1,1,1,1)$\\
466&$>$&$\leq$&$\leq$&$\leq$&$>$&$\leq$&$>$&$>$&$>$&$\leq$&1&$(8^*,2)$\\
467&$>$&$\leq$&$\leq$&$\leq$&$>$&$\leq$&$>$&$>$&$\leq$&$>$&5(ii)&$(2,1,\cdots,1),(3,1,2,1,2,1)$\\
468&$>$&$\leq$&$\leq$&$\leq$&$>$&$\leq$&$>$&$>$&$\leq$&$\leq$&1&$(7^*,3)$\\
469&$>$&$\leq$&$\leq$&$\leq$&$>$&$\leq$&$>$&$\leq$&$>$&$>$&5(ii)&$(2,1,\cdots,1),(3,1,2,2,1,1)$\\
470&$>$&$\leq$&$\leq$&$\leq$&$>$&$\leq$&$>$&$\leq$&$>$&$\leq$&1&$(8^*,2)$\\
471&$>$&$\leq$&$\leq$&$\leq$&$>$&$\leq$&$>$&$\leq$&$\leq$&$>$&5(i)&$(2,1,\cdots,1),(4,1,\cdots,1)$\\
472&$>$&$\leq$&$\leq$&$\leq$&$>$&$\leq$&$>$&$\leq$&$\leq$&$\leq$&2&$(2,1,1,1,1,1,1,1,1)$\\
473&$>$&$\leq$&$\leq$&$\leq$&$>$&$\leq$&$\leq$&$>$&$>$&$>$&5(ii)&$(2,1,\cdots,1),(3,1,3,1,1,1)$ \\
474&$>$&$\leq$&$\leq$&$\leq$&$>$&$\leq$&$\leq$&$>$&$>$&$\leq$&1&$(8^*,2)$ \\
475&$>$&$\leq$&$\leq$&$\leq$&$>$&$\leq$&$\leq$&$>$&$\leq$&$>$&5(i)&$(2,1,\cdots,1),(4,1,\cdots,1)$ \\
476&$>$&$\leq$&$\leq$&$\leq$&$>$&$\leq$&$\leq$&$>$&$\leq$&$\leq$&1&$(7^*,3)$ \\
477&$>$&$\leq$&$\leq$&$\leq$&$>$&$\leq$&$\leq$&$\leq$&$>$&$>$&5(i)&$(2,1,\cdots,1),(4,1,\cdots,1)$ \\
478&$>$&$\leq$&$\leq$&$\leq$&$>$&$\leq$&$\leq$&$\leq$&$>$&$\leq$&1&$(8^*,2)$ \\
479&$>$&$\leq$&$\leq$&$\leq$&$>$&$\leq$&$\leq$&$\leq$&$\leq$&$>$&2&$(2,1,1,1,1,1,1,1,1)$ \\
480&$>$&$\leq$&$\leq$&$\leq$&$>$&$\leq$&$\leq$&$\leq$&$\leq$&$\leq$&2&$(2,1,1,1,1,1,1,1,1)$ \\
\end{tabular}
\newpage
 \begin{tabular}{ccccccccccccc}
 Case & A & B & C & D & E & F & G & H & I &J&Proposition& Inequalities \\ &&&&&&&&&&&&\\
481&$>$&$\leq$&$\leq$&$\leq$&$\leq$&$>$&$>$&$>$&$>$&$>$&-&$(3,1,\cdots,1)$ \\
482&$>$&$\leq$&$\leq$&$\leq$&$\leq$&$>$&$>$&$>$&$>$&$\leq$&1&$(8^*,2)$ \\
483&$>$&$\leq$&$\leq$&$\leq$&$\leq$&$>$&$>$&$>$&$\leq$&$>$&-&$(3,1,1,1,1,2,1),(4,1,1,1,2,1)$ \\
484&$>$&$\leq$&$\leq$&$\leq$&$\leq$&$>$&$>$&$>$&$\leq$&$\leq$&1&$(7^*,3)$ \\
485&$>$&$\leq$&$\leq$&$\leq$&$\leq$&$>$&$>$&$\leq$&$>$&$>$&-&$(3,1,1,1,2,1,1)$,$(4,1,1,2,1,1)$ \\
486&$>$&$\leq$&$\leq$&$\leq$&$\leq$&$>$&$>$&$\leq$&$>$&$\leq$&1&$(8^*,2)$ \\
487&$>$&$\leq$&$\leq$&$\leq$&$\leq$&$>$&$>$&$\leq$&$\leq$&$>$&5(i)&$(2,1,\cdots,1),(4,1,\cdots,1)$\\
488&$>$&$\leq$&$\leq$&$\leq$&$\leq$&$>$&$>$&$\leq$&$\leq$&$\leq$&2&$(2,1,1,1,1,1,1,1,1)$\\
489&$>$&$\leq$&$\leq$&$\leq$&$\leq$&$>$&$\leq$&$>$&$>$&$>$&-&$(3,1,1,2,1,1,1)$,\\
490&$>$&$\leq$&$\leq$&$\leq$&$\leq$&$>$&$\leq$&$>$&$>$&$\leq$&1&$(8^*,2)$\\
491&$>$&$\leq$&$\leq$&$\leq$&$\leq$&$>$&$\leq$&$>$&$\leq$&$>$&5(i)&$(2,1,\cdots,1),(4,1,\cdots,1)$ \\
492&$>$&$\leq$&$\leq$&$\leq$&$\leq$&$>$&$\leq$&$>$&$\leq$&$\leq$&1&$(7^*,3)$ \\
493&$>$&$\leq$&$\leq$&$\leq$&$\leq$&$>$&$\leq$&$\leq$&$>$&$>$&5(i)&$(2,1,\cdots,1),(4,1,\cdots,1)$ \\
494&$>$&$\leq$&$\leq$&$\leq$&$\leq$&$>$&$\leq$&$\leq$&$>$&$\leq$&1&$(8^*,2)$ \\
495&$>$&$\leq$&$\leq$&$\leq$&$\leq$&$>$&$\leq$&$\leq$&$\leq$&$>$&2&$(2,1,1,1,1,1,1,1,1)$ \\
496&$>$&$\leq$&$\leq$&$\leq$&$\leq$&$>$&$\leq$&$\leq$&$\leq$&$\leq$&2&(2,1,1,1,1,1,1,1,1) \\
497&$>$&$\leq$&$\leq$&$\leq$&$\leq$&$\leq$&$>$&$>$&$>$&$>$&-&$(2,1,\cdots,1)$\\
498&$>$&$\leq$&$\leq$&$\leq$&$\leq$&$\leq$&$>$&$>$&$>$&$\leq$&1&$(8^*,2)$ \\
499&$>$&$\leq$&$\leq$&$\leq$&$\leq$&$\leq$&$>$&$>$&$\leq$&$>$&5(i)&$(2,1,\cdots,1),(4,1,\cdots,1)$ \\
500&$>$&$\leq$&$\leq$&$\leq$&$\leq$&$\leq$&$>$&$>$&$\leq$&$\leq$&1&$(7^*,3)$ \\
501&$>$&$\leq$&$\leq$&$\leq$&$\leq$&$\leq$&$>$&$\leq$&$>$&$>$&5(i)&$(2,1,\cdots,1),(4,1,\cdots,1)$ \\
502&$>$&$\leq$&$\leq$&$\leq$&$\leq$&$\leq$&$>$&$\leq$&$>$&$\leq$&1&$(8^*,2)$ \\
503&$>$&$\leq$&$\leq$&$\leq$&$\leq$&$\leq$&$>$&$\leq$&$\leq$&$>$&2&$(2,1,1,1,1,1,1,1,1)$ \\
504&$>$&$\leq$&$\leq$&$\leq$&$\leq$&$\leq$&$>$&$\leq$&$\leq$&$\leq$&2&$(2,1,1,1,1,1,1,1,1)$ \\
505&$>$&$\leq$&$\leq$&$\leq$&$\leq$&$\leq$&$\leq$&$>$&$>$&$>$&5(i)&$(2,1,\cdots,1),(4,1,\cdots,1)$ \\
506&$>$&$\leq$&$\leq$&$\leq$&$\leq$&$\leq$&$\leq$&$>$&$>$&$\leq$&1&$(8^*,2)$ \\
507&$>$&$\leq$&$\leq$&$\leq$&$\leq$&$\leq$&$\leq$&$>$&$\leq$&$>$&2&$(2,1,1,1,1,1,1,1,1)$ \\
508&$>$&$\leq$&$\leq$&$\leq$&$\leq$&$\leq$&$\leq$&$>$&$\leq$&$\leq$&1&$(7^*,3)$ \\
509&$>$&$\leq$&$\leq$&$\leq$&$\leq$&$\leq$&$\leq$&$\leq$&$>$&$>$&2&$(2,1,1,1,1,1,1,1,1)$ \\
510&$>$&$\leq$&$\leq$&$\leq$&$\leq$&$\leq$&$\leq$&$\leq$&$>$&$\leq$&1&$(8^*,2)$ \\
511&$>$&$\leq$&$\leq$&$\leq$&$\leq$&$\leq$&$\leq$&$\leq$&$\leq$&$>$&2&$(2,1,1,1,1,1,1,1,1)$ \\
512&$>$&$\leq$&$\leq$&$\leq$&$\leq$&$\leq$&$\leq$&$\leq$&$\leq$&$\leq$&2&$(2,1,1,1,1,1,1,1,1)$ \\
\end{tabular}}

\section{Easy Cases}\label{Easy}
\begin{prop} Cases where $I > 1,~J \leq 1$ or where $H>1, ~I\leq 1, ~J \leq 1$ do not arise.\end{prop}
\noindent Proof is similar to that of Propositions 1 and 2 of \cite{HRS8}.\vspace{1mm}\\
 This settles $128+64$ cases.\hfill$\Box$
 \begin{prop} Cases in which $B\leq1$ and at most two out of $C,D,E,F,G,H,I,J$ are greater than 1 do not arise.\end{prop}
\noindent Proof is similar to that of Proposition 3(i) of  \cite{HRS8}.\vspace{1mm}\\
 This settles another 23 cases.\hfill$\Box$\\
\begin{lemma}\label{lemma 6} Let $X_{1},\cdots,X_{10}$ be positive real numbers, each $<2.2636302$ and satisfying
\begin{equation}\label{1.4.1}X_{1}>1,~~~X_{1}X_{2}\cdots X_{10}=1~~~{\rm and}~~~x_i=|X_i-1|.\end{equation}
Then the following hold :\vspace{2mm}\\
(i)~~~{If $X_{i}>1$ for $3\leq i\leq10$, then we have}\vspace{1mm}\\
\indent~~ $\mathfrak{S}_1=4X_1-\frac{2X_{1}^2}{X_2}+X_3+\cdots+X_{10}\leq10$.\vspace{2mm}\\
(ii)~~~{If $X_{i}>1$ for $i=3,5,6,7,8,9,10$, then we have}\vspace{1mm}\\
\indent~~  $\mathfrak{S}_2=4X_1-\frac{2X_{1}^2}{X_2}+4X_3-\frac{2X_{3}^2}{X_4}+X_5+\cdots+X_{10}\leq10$.\vspace{2mm}\\
(iii)~~~{If $X_{i}>1$ for $i=3,5,7,8,9,10$, then we have}\vspace{1mm}\\
\indent~~  $\mathfrak{S}_3=4X_1-\frac{2X_{1}^2}{X_2}+4X_3-\frac{2X_{3}^2}{X_4}+4X_5-\frac{2X_{5}^2}{X_6}+X_7+X_8+X_9+X_{10}\leq10$.\vspace{2mm}\\
(iv)~~~{If $X_{i}>1$ for $i=3,5,7,9,10$, then we have}\vspace{1mm}\\
\indent~~  $\mathfrak{S}_4=4X_1-\frac{2X_{1}^2}{X_2}+4X_3-\frac{2X_{3}^2}{X_4}+4X_5-\frac{2X_{5}^2}{X_6}+4X_7-\frac{2X_{7}^2}{X_8}+X_9+X_{10}\leq10$.\vspace{2mm}\\
(v)~~~{If $X_{i}>1$ for $i=4,5,6,7,8,9,10$, then we have}\vspace{1mm}\\
\indent~~  $\mathfrak{S}_5=4X_1-\frac{X_{1}^3}{X_{2}X_{3}}+X_4+X_{5}+X_{6}+X_7+X_8+X_9+X_{10}\leq10$.\vspace{2mm}\\
(vi)~~~{If $X_{i}>1$ for $i=4,7,8,9,10$ and $X_7\leq X_1$, $X_8\leq X_1$, $X_9\leq X_1$, $X_{10}\leq X_1$}\vspace{1mm}\\
\indent~~ then we have\vspace{1mm}\\
\indent~~  $\mathfrak{S}_6=4X_1-\frac{X_{1}^3}{X_{2}X_{3}}+4X_4-\frac{X_{4}^3}{X_{5}X_{6}}+X_7+X_8+X_9+X_{10}\leq10$.\vspace{2mm}\\
(vii)~~~{If $X_{i}>1$ for $i=4,7,10$ and $X_{10}\leq X_1$, then we have}\vspace{1mm}\\
\indent~~  $\mathfrak{S}_7=4X_1-\frac{X_{1}^3}{X_{2}X_{3}}+4X_4-\frac{X_{4}^3}{X_{5}X_{6}}+4X_7-\frac{X_{7}^3}{X_{8}X_{9}}+X_{10}\leq10$.\vspace{2mm}\\
(viii)~~~{If $X_{i}>1$ for $i=4,6,7,8,9,10$ and $X_i\leq X_{1}X_{4}$ for $i=6,7,8,9,10$,\\\indent~~ $X_5\leq X_4$, }
  then we have\vspace{1mm}\\
\indent~~  $\mathfrak{S}_8=4X_1-\frac{X_{1}^3}{X_{2}X_{3}}+4X_4-\frac{2X_{4}^2}{X_{5}}+X_6+X_7+X_8+X_9+X_{10}\leq10$.\vspace{2mm}\\
(ix)~~~{If $X_{i}>1$ for $i=3,5,8,9,10$ and $X_2\leq X_1$, $X_4\leq X_3$, $X_i\leq X_{1}X_{5}$, \\\indent~~for $i=8,9,10$,} then we have\vspace{1mm}\\
\indent~~  $\mathfrak{S}_9=4X_1-\frac{2X_{1}^2}{X_2}+4X_3-\frac{2X_{3}^2}{X_4}+4X_5-\frac{X_{5}^3}{X_{6}X_{7}}+X_8+X_9+X_{10}\leq10$.\vspace{2mm}\\
(x)~~~{If $X_{i}>1$ for $i=3,5,7,10$ and $X_2\leq X_1$, $X_4\leq X_3$, $X_6\leq X_5$, $X_{10}\leq X_{1}X_{7}$,\\\indent~~ then we have}\vspace{1mm}\\
\indent~~  $\mathfrak{S}_{10}=4X_1-\frac{2X_{1}^2}{X_2}+4X_3-\frac{2X_{3}^2}{X_4}+4X_5-\frac{2X_{5}^2}{X_6}+4X_7-\frac{X_{7}^3}{X_{8}X_{9}}+X_{10}\leq10$.\vspace{2mm}\\
(xi)~~~{If $X_{i}>1$ for $i=4,7,9,10$, $X_8\leq X_7$, $X_9\leq X_1X_{4}X_{7}$, $X_{10}\leq X_1X_{4}X_{7}$,\\\indent~~ then we have}\vspace{1mm}\\
\indent~~  $\mathfrak{S}_{11}=4X_1-\frac{X_{1}^3}{X_{2}X_{3}}+4X_4-\frac{X_{4}^3}{X_{5}X_{6}}+4X_7-\frac{2X_{7}^2}{X_8}+X_9+X_{10}\leq10$.\end{lemma}
Proof is similar to that of Lemma 5 of  \cite{HRS8}, so omitted.\hfill$\Box$\\
{\noindent \bf Remark.}  The above lemma can be generalized for arbitrary $n$.
\begin{prop}\label{Prop 3} The following cases do not arise.\\
$\begin{array}{ll}{\rm (i)}&(3),(5),(9),(17),(33),(65),(129),(257)\\
{\rm (ii)}&(11), (19), (21), (35), (37), (41), (67), (69), (73), (81), (131), (133), (137), (145),\\& (161), (259), (261), (265), (273), (289), (321).\\
{\rm (iii)}&(43), (75), (83), (85), (139), (147), (149), (163), (165), (169), (267), (275), (277),\\& (291), (293), (297), (323), (325), (329), (337). \\
{\rm (iv)}&(171), (299), (331), (339), (341).\\
{\rm (v)}&(7), (13), (25), (49), (97), (193), (385).\\
{\rm (vi)}&(135), (141), (153), (177), (391), (397), (409), (433).\\
\end{array} $
$\begin{array}{ll}
{\rm (vii)}&(55), (183), (439).\\
{\rm (viii)}&(263), (269), (281), (305), (353), (387), (389), (393), (401), (417).\\
{\rm (ix)}&(23), (71), (77), (89), (99), (101), (105), (155), (167), (173), (179), (181), (279), (283),\\& (295), (301), (307), (309), (327), (345), (355), (357), (361), (395), (403), (405),\\& (419), (421), (425).\\
{\rm (x)}&(87), (91), (107), (343), (347), (363), (427).\\
{\rm (xi)}&(27), (39), (45), (51), (53), (103), (109), (151), (311), (359), (365), (407), (411),\\& (423), (429), (435), (437).\end{array} $
\end{prop}
Each part of Proposition 3 follows from corresponding part of Lemma 6, after selecting suitable inequality. The inequality used for each case has been mentioned in Table I.\hfill$\Box$
\begin{lemma} Let $X_{1},\cdots,X_{10}$ be positive real numbers, each $< 2.2636302$, satisfying $\eqref{1.4.1}$. Let
$$ \gamma={\displaystyle\sum_{4 \leq i \leq 10\atop{X_{i} \leq 1}}}x_{i}~~~{\rm and} ~~~~\delta={\displaystyle\sum_{4 \leq i
 \leq 10\atop{X_{i}>1}}}x_{i}$$
Suppose that $\gamma\leq x_{1}\leq 0.5$,
then $$\mathfrak{S}_{12} = 4X_1-X_{1}^{4}X_{4}\ldots X_{10}+X_{4}+\ldots +X_{10}\leq 10.$$\end{lemma}
Proof is similar to that of Lemma 6(i) of  \cite{HRS8}, so omitted.\hfill$\Box$\\
\begin{prop}
  The Cases 216, 224, 240 do not arise.\end{prop}
{\noindent Proof.} These cases are dealt by using the above lemma. The suitable inequalities used for each case, have been mentioned in the Table I.\hfill$\Box$


\section{Difficult Cases}\label{Difficult}
In this section we consider 151 cases, which need more intricate analysis of available inequalities. These cases have been solved  using the Optimization tool of the software Mathemetica(non-linear global optimization) for the optimization of the functions in 10 variables $A, B, C, D, E, F, G, H, I, J$ under a number of constraints.
\par
\begin{figure}[h!]
\includegraphics[width=155mm,height=155mm]{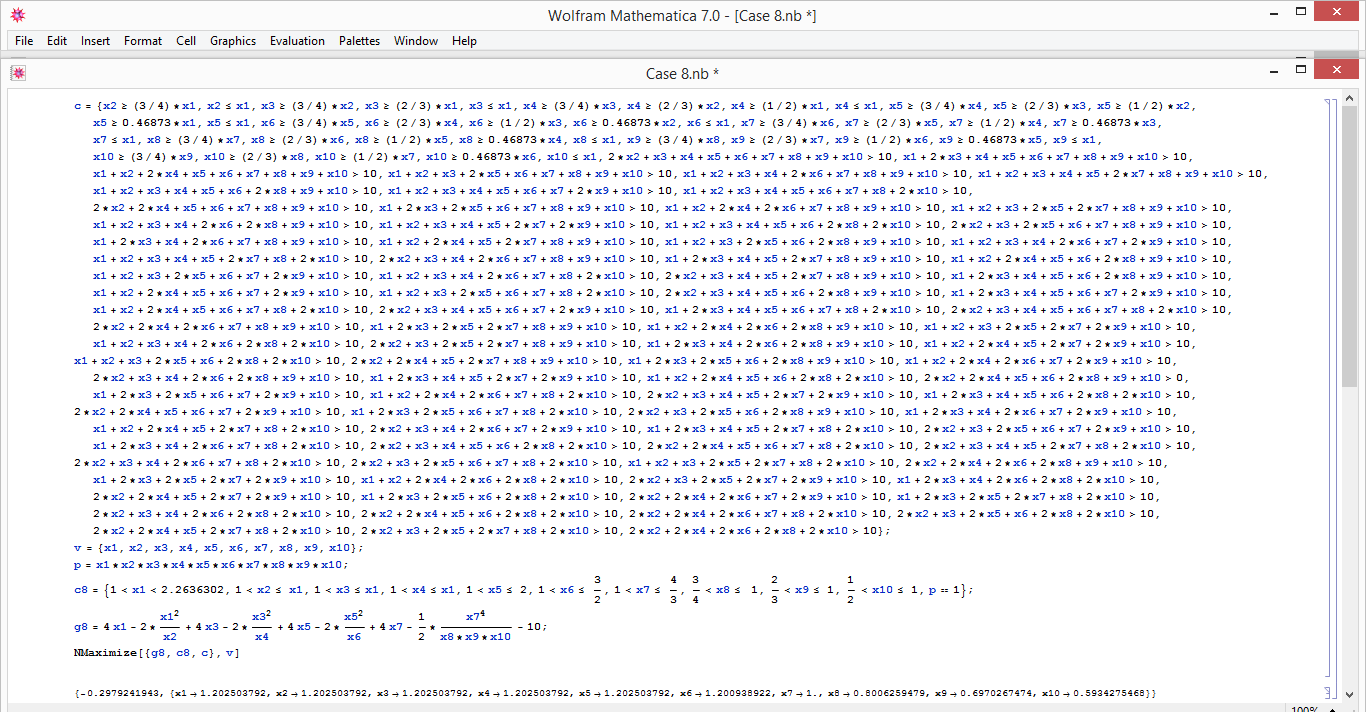}
\caption{Mathematica Model for $\phi_3$ (Case 8)}
\label{Fig 1}
\end{figure}

\begin{figure}[h!]
\includegraphics[width=155mm,height=155mm]{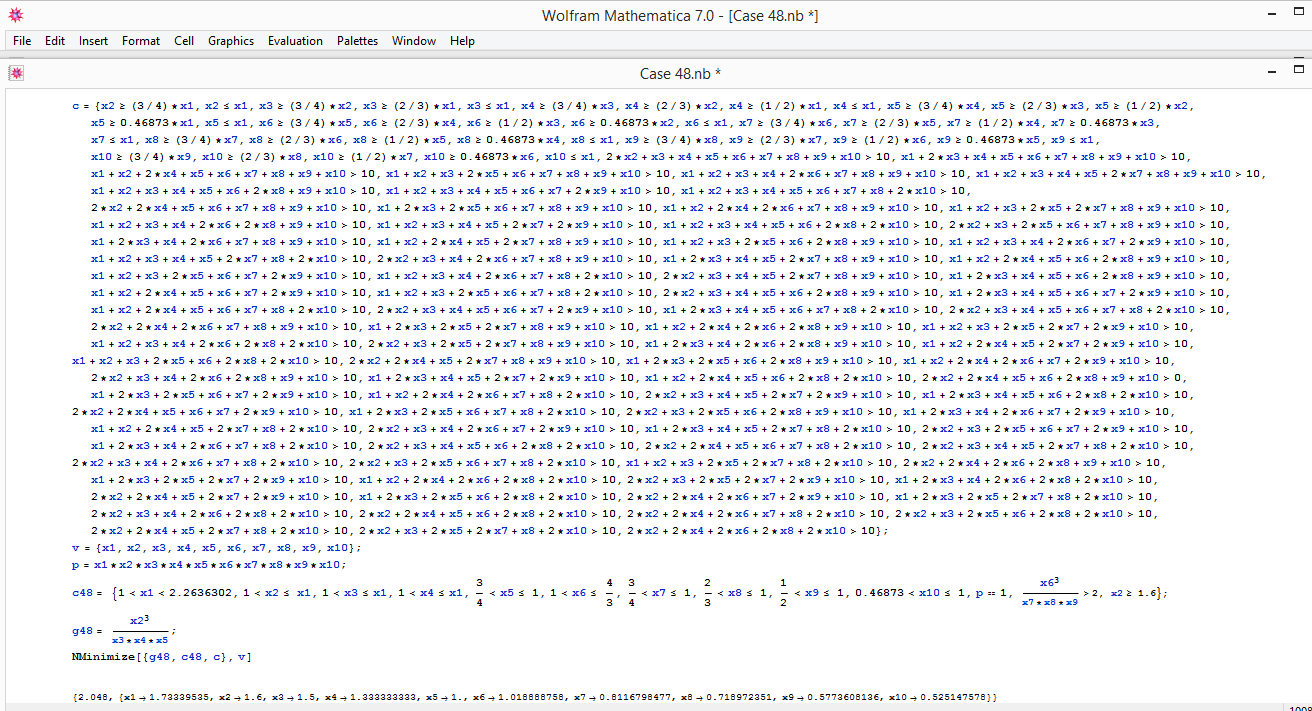}
\caption{Mathematica Model for $\psi_2$ (Case 48)}
\label{Fig 2}
\end{figure}

 NMinimize and NMaximize implement several algorithms for finding constrained global optima. The methods are flexible enough to cope with functions that are not differentiable or continuous and are not easily trapped by local optima. The constraints to NMinimize and NMaximize may be either a list or a logical combination of equalities, inequalities, and domain specifications. Equalities and inequalities may be nonlinear. The settings for Accuracy Goal and Precision Goal specify the number of digits to seek in both the value of the position of the maximum, and the value of the function at the maximum. NMaximize continues until either of the goals specified by AccuracyGoal or Precision Goal is achieved. The default settings for Accuracy Goal and Precision Goal are Working Precision/2. Method is an option for various algorithm-intensive functions that specifies what internal methods they should use. With the default setting Method$->$Automatic Mathematica will automatically try to pick the best method for a particular computation.\vspace{2mm}\\
The common and obvious constraints, which arise from Lemmas 3 and 4, for all cases are :\vspace{2mm}\\
$\begin{array}{lllll}
1<A<2.2636302,&B\geq(3/4)\times A,&B\leq A,&&\\
C\geq(3/4)\times B,&C\geq(2/3)\times A,&C\leq A,&&\\
D\geq(3/4)\times C,&D\geq(2/3)\times B,&D\geq(1/2)\times A,&D\leq A,&\\
E\geq(3/4)\times D,&E\geq(2/3)\times C,&E\geq(1/2)\times B,& E\geq0.46873\times A,& E\leq A,\\
F\geq(3/4)\times E,&F\geq(2/3)\times D,&F\geq(1/2)\times C,&F\geq0.46873\times B,&F\leq A,\\
G\geq(3/4)\times F,&G\geq(2/3)\times E,&G\geq(1/2)\times D,&G\geq0.46873\times C,&G\leq A,\\
H\geq(3/4)\times G,&H\geq(2/3)\times F,&H\geq(1/2)\times E,&H\geq0.46873\times D,&H\leq A,\\
I\geq(3/4)\times H,&I\geq(2/3)\times G,&I\geq(1/2)\times F ,&I\geq0.46873\times E,&I\leq A,\\
J\geq(3/4)\times I,&J\geq(2/3)\times H,&J\geq(1/2)\times G,&J\geq0.46873\times F,&J\leq A.\\
\end{array}$\vspace{2mm}\\
We also use the following 88 constraints which arise from all the possible \emph{weak} inequalities
corresponding to the partitions of 10 with summands equal to 1 and 2:\vspace{2mm}\\
$\begin{array}{ll}
2B+C+D+E+F+G+H+I+J>10,& A+2C+D+E+F+G+H+I+J>10,\\
A+B+2D+E+F+G+H+I+J>10,&A+B+C+2E+F+G+H+I+J>10,\\
 A+B+C+D+2F+G+H+I+J>10,& A+B+C+D+E+2G+H+I+J>10,\\
A+B+C+D+E+F+2H+I+J>10,& A+B+C+D+E+F+G+2I+J>10,\\
 A+B+C+D+E+F+G+H+2J>10,&\\
2B+2D+E+F+G+H+I+J>10,
&A+2C+2E+F+G+H+I+J>10,
\\A+B+2D+2F+G+H+I+J>10,
&A+B+C+2E+2G+H+I+J>10,
\\ A+B+C+D+2F+2H+I+J>10,
 &A+B+C+D+E+2G+2I+J>10,
\\A+B+C+D+E+F+2H+2J>10,
&2B+C+2E+F+G+H+I+J>10,
\\A+2C+D+2F+G+H+I+J>10,
&A+B+2D+E+2G+H+I+J>10,
\\A+B+C+2E+F+2H+I+J>10,
&A+B+C+D+2F+G+2I+J>10,
  \\A+B+C+D+E+2G+H+2J>10,
&2B+C+D+2F+G+H+I+J>10,
\\A+2C+D+E+2G+H+I+J>10,
 &A+B+2D+E+F+2H+I+J>10,
 \\A+B+C+2E+F+G+2I+J>10,
&A+B+C+D+2F+G+H+2J>10, 
\\2B+C+D+E+2G+H+I+J>10,
 &A+2C+D+E+F+2H+I+J>10,
 \\A+B+2D+E+F+G+2I+J>10,
&A+B+C+2E+F+G+H+2J>10,
\\2B+C+D+E+F+2H+I+J>10,
&A+2C+D+E+F+G+2I+J>10,
\\A+B+2D+E+F+G+H+2J>10,
&2B+C+D+E+F+G+2I+J>10,
 \\A+2C+D+E+F+G+H+2J>10, 
 &2B+C+D+E+F+G+H+2J>10, 
\\2B+2D+2F+G+H+I+J>10,
 &A+2C+2E+2G+H+I+J>10,
 \\A+B+2D+2F+2H+I+J>10,
&A+B+C+2E+2G+2I+J>10,
\\A+B+C+D+2F+2H+2J>10,
&2B+C+2E+2G+H+I+J>10,
\\A+2C+D+2F+2H+I+J>10,
&A+B+2D+E+2G+2I+J>10,
\\A+B+C+2E+F+2H+2J>10, 
&2B+2D+E+2G+H+I+J>10,
\\A+2C+2E+F+2H+I+J>10,
&A+B+2D+2F+G+2I+J>10,
\\2B+C+D+2F+2H+I+J>10,
&A+2C+D+E+2G+2I+J>10,
\\A+B+2D+E+F+2H+2J>10,
&2B+2D+E+F+2H+I+J>10,
\\A+2C+2E+F+G+2I+J>10,
&A+B+2D+2F+G+H+2J>10, 
\\2B+C+D+E+2G+2I+J>10,
&A+2C+D+E+F+2H+2J>10, 
\\2B+2D+E+F+G+2I+J>10,
&A+2C+2E+F+G+H+2J>10, 
\\2B+C+2E+F+2H+I+J>10,
&A+2C+D+2F+G+2I+J>10,
\\A+B+2D+E+2G+H+2J>10,
&2B+C+D+2F+G+2I+J>10,
\\A+2C+D+E+2G+H+2J>10, 
&2B+C+2E+F+G+2I+J>10,
\\A+2C+D+2F+G+H+2J>10,
&2B+C+D+E+F+2H+2J>10, 
\\2B+2D+E+F+G+H+2J>10, 
&2B+C+D+E+2G+H+2J>10,
  \end{array}$
 $\begin{array}{ll}
\\2B+C+D+2F+G+H+2J>10,
&2B+C+2E+F+G+H+2J>10, 
\\A+B+C+2E+2G+H+2J>10,
&2B+2D+2F+2H+I+J>10,
\\A+2C+2E+2G+2I+J>10,
&A+B+2D+2F+2H+2J>10, 
\\2B+C+2E+2G+2I+J>10,
 &A+2C+D+2F+2H+2J>10, 
\\2B+2D+E+2G+2I+J>10,
 &A+2C+2E+F+2H+2J>10, 
\\2B+2D+2F+G+2I+J>10,
&A+2C+2E+2G+H+2J>10,
\\2B+C+D+2F+2H+2J>10,
&2B+2D+E+F+2H+2J>10, 
\\2B+2D+2F+G+H+2J>10, 
&2B+C+2E+F+2H+2J>10, 
\\2B+2D+E+2G+H+2J>10, 
&2B+C+2E+2G+H+2J>10,
\\2B+2D+2F+2H+2J>10. &
\end{array}$\vspace{2mm}\\
Some specific constraints due to the numerical bounds on the variables arise in the
individual case.
In order to write all these constraints in Mathematica Notebook, we call $A=x1, B=x2, \cdots, J=x10.$ Let $\phi_l$ and $\psi_l$ for $l\geq1$,
be non-linear functions in 10 variables $A, B, \cdots, J$, listed below. The functions $\phi_l$'s are associated with various conditional inequalities and the same have been written against the respective function. In each of the cases, we find using Mathematica that Global maximum of $\phi_l$ is always
less than 10 and the Global minimum of  $\psi_l$ is always greater than 2 under all the
above mentioned constraints and some specific constraints relating to the numerical
bounds on the variables, which consequently contradicts to the corresponding inequality (as mentioned in Section \ref{Plan}).\\\hyperref[Fig 1]{Figure 1}  shows the Mathematica Model of $\phi_3$ used in the Case 8 and \hyperref[Fig 2]{Figure 2} shows Mathematica Model of $\psi_2$ used in Case 48.\vspace{2mm}\\
$\psi_1\label{s1}=\frac{A^3}{BCD}$;~$\psi_2\label{s2}=\frac{B^3}{CDE}$;~$\psi_3\label{s3}=\frac{C^3}{DEF}$;~$\psi_4\label{s4}=\frac{D^3}{EFG}$;
~$\psi_5\label{s5}=\frac{E^3}{FGH}$;~$\psi_6\label{s6}=\frac{F^3}{GHI}$;~$\psi_7\label{s7}=\frac{G^3}{HIJ}$;
$\begin{array}{ll}
\phi_1\label{k1}=4A-\frac{2A^2}{B}+C+D+4E-\frac{E^3}{FG}+H+I+J ;& (2,1,1,3,1,1,1)\vspace{1mm}\\
\phi_2\label{k2}=4A-\frac{2A^2}{B}+C+D+E+4F-\frac{F^3}{GH}+I+J ;& (2,1,1,1,3,1,1)\vspace{1mm}\\
\phi_3\label{h1}=4A-\frac{2A^2}{B}+4C-\frac{2C^2}{D}+4E-\frac{2E^2}{F}+4G-\frac{1}{2}\frac{G^4}{HIJ} ;& (2,2,2,4)\vspace{1mm}\\
\phi_4\label{h2}=4A-\frac{A^3}{BC}+D+E+F+4G-\frac{G^3}{HI}+J ;& (3,1,1,1,3,1)\vspace{1mm}\\
\phi_5\label{h3}=4A-\frac{2A^2}{B}+4C-\frac{2C^2}{D}+E+4F-\frac{1}{2}\frac{F^4}{GHI}+J ; &(2,2,1,4,1)\vspace{1mm}\\
\phi_6\label{h4}=4A-\frac{A^3}{BC}+D+E+4F-\frac{F^3}{GH}+I+J ;&(3,1,1,3,1,1)\vspace{1mm}\\
\phi_7\label{h5}=4A-\frac{2A^2}{B}+4C-\frac{2C^2}{D}+4E-\frac{1}{2}\frac{E^4}{FGH}+I+J ; &(2,2,4,1,1)\vspace{1mm}\\
\phi_8\label{h6}=4A-\frac{A^3}{BC}+D+4E-\frac{E^3}{FG}+H+I+J ;& (3,1,3,1,1,1)\vspace{1mm}\\
\phi_9\label{h7}=A+4B-\frac{2B^2}{C}+4D-\frac{2D^2}{E}+4F-\frac{1}{2}\frac{F^4}{GHI}+J ; &(1,2,2,4,1)\vspace{1mm}\\
\phi_{10}\label{h8}=A+4B-\frac{1}{2}\frac{B^4}{CDE}+4F-\frac{1}{2}\frac{F^4}{GHI}+J ;& (1,4,4,1)\vspace{1mm}\\
\phi_{11}\label{h9}=A+4B-\frac{2B^2}{C}+4D-\frac{2D^2}{E}+4F-\frac{2F^2}{G}+4H-\frac{2H^2}{I}+J ; &(1,2,2,2,2,1)\vspace{1mm}\\
\phi_{12}\label{h10}=4A-\frac{1}{2}\frac{A^4}{BCD}+4E-\frac{2E^2}{F}+4G-\frac{2G^2}{H}+I+J ;& (4,2,2,1,1)\vspace{1mm}\\
\end{array}$
\newpage
$\begin{array}{ll}
\phi_{13}\label{h11}=A+4B-\frac{1}{2}\frac{B^4}{CDE}+4F-\frac{2F^2}{G}+4H-\frac{2H^2}{I}+J ; &(1,4,2,2,1)\vspace{1mm}\\
\phi_{14}\label{h12}=4A-\frac{2A^2}{B}+4C-\frac{2C^2}{D}+4E-\frac{2E^2}{F}+4G-\frac{2G^2}{H}+I+J ;& (2,2,2,2,1,1)\vspace{1mm}\\
\phi_{15}\label{h13}=4A-\frac{2A^{2}}{B}+4C-\frac{1}{2}\frac{C^4}{DEF}+4G-\frac{2G^2}{H}+I+J ;& (2,4,2,1,1)\vspace{1mm}\\
\phi_{16}\label{h14}=A+4B-\frac{1}{2}\frac{B^4}{CDE}+4F-\frac{2F^2}{G}+H+I+J ;& (1,4,2,1,1,1)\vspace{1mm}\\
\phi_{17}\label{h15}=4A-\frac{1}{2}\frac{A^4}{BCD}+4E-\frac{2E^2}{F}+G+H+I+J ;& (4,2,1,1,1,1)\vspace{1mm}\\
\phi_{18}\label{h16}=4A-\frac{2A^2}{B}+4C-\frac{2C^2}{D}+4E-\frac{2E^2}{F}+G+H+I+J ;& (2,2,2,1,1,1,1)\vspace{1mm}\\
\phi_{19}\label{h17}=4A-\frac{2A^2}{B}+4C-\frac{2C^2}{D}+4E-\frac{2E^2}{F}+G+4H-\frac{2H^2}{I}+J ;& (2,2,2,1,2,1)\vspace{1mm}\\
\phi_{20}\label{h18}=4A-\frac{2A^2}{B}+4C-\frac{1}{2}\frac{C^4}{DEF}+G+4H-\frac{2H^2}{I}+J ;& (2,4,1,2,1)\vspace{1mm}\\
\phi_{21}\label{h19}=A+4B-\frac{1}{2}\frac{B^4}{CDE}+F+G+H+I+J ;& (1,4,1,1,1,1)\vspace{1mm}\\
\phi_{22}\label{h20}=A+4B-\frac{2B^2}{C}+4D-\frac{D^3}{EF}+G+H+I+J ;& (1,2,3,1,1,1,1)\vspace{1mm}\\
\phi_{23}\label{h21}=4A-\frac{A^3}{BC}+D+E+F+G+H+I+J ;& (3,1,\cdots,1)\vspace{1mm}\\
\phi_{24}\label{h22}=4A-\frac{1}{2}\frac{A^4}{BCD}+E+F+G+H+I+J ;& (4,1,\cdots,1)\vspace{1mm}\\
\phi_{25}\label{h23}=4A-\frac{2A^2}{B}+4C-\frac{2C^2}{D}+4E-\frac{1}{2}\frac{E^4}{FGH}+4I-\frac{2I^2}{J} ;& (2,2,4,2)\vspace{1mm}\\
\phi_{26}\label{h24}=4A-\frac{2A^2}{B}+C+4D-\frac{1}{2}\frac{D^4}{EFG}+4H-\frac{2H^2}{I}+J ;& (2,1,4,2,1)\vspace{1mm}\\
\phi_{27}\label{h25}=4A-\frac{2A^2}{B}+4C-\frac{1}{2}\frac{C^4}{DEF}+4G-\frac{1}{2}\frac{G^4}{HIJ} ;& (2,4,4)\vspace{1mm}\\
\phi_{28}\label{h26}=4A-\frac{2A^2}{B}+4C-\frac{1}{2}\frac{C^4}{DEF}+4G-\frac{2G^2}{H}+4I-\frac{2I^2}{J} ;& (2,4,2,2)\vspace{1mm}\\
\phi_{29}\label{h27}=A+4B-\frac{2B^2}{C}+4D-\frac{1}{2}\frac{D^4}{EFG}+4H-\frac{2H^2}{I}+J ;& (1,2,4,2,1)\vspace{1mm}\\
\phi_{30}\label{h28}=4A-\frac{2A^2}{B}+4C-\frac{1}{2}\frac{C^4}{DEF}+4G-\frac{2G^2}{H}+I+J ;& (2,4,2,1,1)\vspace{1mm}\\
\phi_{31}\label{h29}=4A-\frac{1}{2}\frac{A^4}{BCD}+4E-\frac{1}{2}\frac{E^4}{FGH}+4I-\frac{2I^2}{J} ;& (4,4,2)\vspace{1mm}\\
\phi_{32}\label{h30}=4A-\frac{1}{2}\frac{A^4}{BCD}+4E-\frac{2E^2}{F}+4G-\frac{2G^2}{H}+4I-\frac{2I^2}{J} ;& (4,2,2,2)\vspace{1mm}\\
\phi_{33}\label{h31}=A+4B-\frac{2B^2}{C}+4D-\frac{1}{2}\frac{D^4}{EFG}+H+I+J ;& (1,2,4,1,1,1)\vspace{1mm}\\
\phi_{34}\label{h32}=4A-\frac{1}{2}\frac{A^4}{BCD}+4E-\frac{1}{2}\frac{E^4}{FGH}+I+J ;& (4,4,1,1)\vspace{1mm}\\
\phi_{35}\label{h33}=4A-\frac{2A^2}{B}+4C-\frac{2C^2}{D}+4E-\frac{1}{2}\frac{E^4}{FGH}+I+J ;& (2,2,4,1,1)\vspace{1mm}\\
\phi_{36}\label{h34}=4A-\frac{A^3}{BC}+4D-\frac{1}{2}\frac{D^4}{EFG}+4H-\frac{2H^2}{I}+J ;& (3,4,2,1)\vspace{1mm}\\
\phi_{37}\label{h35}=A+4B-\frac{2B^2}{C}+4D-\frac{2D^2}{E}+4F-\frac{2F^2}{G}+H+I+J ;& (1,2,2,2,1,1,1).
\end{array}$

\begin{lemma} Cases in which $B\leq1$ and \\
(i) any three out of $C,D,E,F,G,H,I,J$ are greater than 1 and $A<1.196$ or \\
(ii) any four out of $C,D,E,F,G,H,I,J$ are greater than 1 and $A<1.096$,\\
 do not arise.\end{lemma}
 Proof is similar to that of Proposition 3 of  \cite{HRS8}, so omitted.\hfill$\Box$\\
\begin{prop} The following cases do not arise.\\
$\begin{array}{ll}{\rm (i)}& 288, 304, 312, 319, 336, 344, 351, 360, 367, 375, 379, 381, 400, 408, 415, 431, 456, 463, \\&471, 475, 477, 487, 491, 493, 499, 501, 505.\\
{\rm (ii)}& 272, 280, 287, 296, 303, 315, 317, 328, 335, 349, 371, 373, 377, 392, 399, 413, 441, 455,\\& 459, 461, 467, 469, 473.\end{array}$\end{prop}
{\noindent\bf Proof.} We apply Lemma 8(i) for the cases in part (i) and Lemma 8(ii) for the cases in part (ii). For sake of convenience we illustrate the Cases 288 and 272 here. Other cases are similar and we  have mentioned the inequalities used against each  in Table I.\\
Case 288: $A>1, B\leq1, C>1, D>1, E>1, F\leq1, G\leq1, H\leq1, I\leq1, J\leq1 $\\
Using the Lemma 8(i) we get $A\geq1.196$. Now using Optimization Tool of Mathematica we find that $max~\hyperref[k1]{\phi_1}\leq10$, which gives contradiction to the corresponding inequality (2,1,1,3,1,1,1).\\
Case 272: $A>1, B\leq1, C>1, D>1, E>1, F>1, G\leq1, H\leq1, I\leq1, J\leq1 $\\
Using the Lemma 8(ii) we get $A\geq1.096$. Now using Optimization Tool of Mathematica we find that $max~\hyperref[k2]{\phi_2}\leq10$, which gives contradiction to the corresponding inequality (2,1,1,1,3,1,1).
\hfill$\Box$\\
Proposition 5 settles 50 difficult cases. Out of remaining 101 difficult  cases, we discuss here 11 cases only in Proposition 6. The other  90 cases are similar and the main  inequalities used to find contradiction are listed in Table 1 against each case.
\noindent
\noindent
\begin{prop}
  The following cases do not arise :
\end{prop}
{\noindent \bf Case 8} $A>1, B>1, C>1, D>1, E>1, F>1, G>1, H\leq1, I\leq1, J\leq1$.\\
Proof.
First suppose $G^3/HIJ\geq2$, then the inequality (2,2,2,4) holds, but $max~\hyperref[h1]{\phi_3}
\leq10$. So $G^3/HIJ<2$, then $max~\hyperref[h2]{\phi_4}\leq10$, which contradicts to the inequality (3,1,1,1,3,1).\vspace{3mm}\\
{\noindent \bf Case 15} $A>1, B>1, C>1, D>1, E>1, F>1, G\leq1, H\leq1, I\leq1, J>1$.\\
Proof.
Suppose $F^3/GHI\geq2$, then the inequality (2,2,1,4,1) holds, but $max~\hyperref[h3]{\phi_5}\leq10$.
So $F^3/GHI<2$, then $max~\hyperref[h4]{\phi_6}\leq10$, which gives contradiction to (3,1,1,3,1,1).\vspace{1mm}\\
{\noindent \bf Case 29} $A>1, B>1, C>1, D>1, E>1, F\leq1, G\leq1, H\leq1, I>1, J>1$.\\
Proof.
Suppose $E^3/FGH\geq2$, then (2,2,4,1,1) holds, but $max~\hyperref[h5]{\phi_7}\leq10$.
So $E^3/FGH<2$, then $max~\hyperref[h6]{\phi_8}\leq10$, which gives contradiction to (3,1,3,1,1,1).\vspace{1mm}\\
{\noindent \bf Case 48} $A>1, B>1, C>1, D>1, E\leq1, F>1, G\leq1, H\leq1, I\leq1, J\leq1$.\\
Proof. Suppose $F^3/GHI>2$. If $B<1.6$, then $max~\hyperref[h7]{\phi_9}\leq10$, contradicts to (1,2,2,4,1) and if $B\geq1.6$, then $ min~\hyperref[s2]{\psi_2}>2$, so the inequality (1,4,4,1) holds, but $max~\hyperref[h8]{\phi_{10}}\leq10$.
Thus we must have $F^3/GHI<2$, but then $max~\hyperref[h9]{\phi_{11}}\leq10$  which is contradiction to the inequality (1,2,2,2,2,1).\vspace{1mm}\\
{\noindent \bf Case 95} $A>1, B>1, C>1, D\leq1, E>1, F\leq1, G\leq1, H\leq1, I\leq1, J>1$.\\
Proof. Claim (i)  $A^3/BCD<2$\\
If $A^3/BCD\geq2$, then (4,2,2,1,1) holds, but  $max~\hyperref[h10]{\phi_{12}}\leq 10$. \\
Claim (ii) $B^3/CDE<2$\\
If $B^3/CDE\geq2$, then the inequality (1,4,2,2,1) holds, but $\hyperref[h11]{\phi_{13}}\leq10$.\\
Finally $max~\hyperref[h12]{\phi_{14}}\leq10$, which gives contradiction to the inequality (2,2,2,2,1,1).\vspace{1mm}\\
{\noindent \bf Case 121} $A>1, B>1, C>1, D\leq1, E\leq1, F\leq1, G\leq1, H>1, I>1, J>1$.\\
Proof: Claim (i) $C<1.05$ \\
If $C\geq1.05$, then $max~\hyperref[s3]{\psi_3}>2$ and so the inequality (2,4,2,1,1) holds, but $max~\hyperref[h13]{\phi_{15}}\leq10$.\\
Claim(ii) $B\leq 1.24$\\
If $B\geq 1.24$, then $min~\hyperref[s2]{\psi_2}>2$, so the inequality (1,4,2,1,1,1) holds, but $max~\hyperref[h14]{\phi_{16}}\leq10$.\\
Claim(iii) $\frac{A^3}{BCD}<2$\\
If $\frac{A^3}{BCD}\geq 2$, then (4,2,1,1,1,1) holds, but $max~\hyperref[h15]{\phi_{17}}\leq10$.\\
Finally the inequality (2,2,2,1,1,1,1) holds, but $max~\hyperref[h16]{\phi_{18}}\leq10$\vspace{1mm}\\
{\noindent \bf Case 123}  $ A > 1, B > 1, C > 1, D \leq 1, E \leq 1, F \leq 1, G \leq 1, H > 1, I \leq 1, J > 1$.\\
Proof: If $C\leq1.08$, then $max~\hyperref[h17]{\phi_{19}}\leq10$, so the inequality (2,2,2,1,2,1) contradicts  and
if $C>1.08$, then $min~\hyperref[s3]{\psi_3}>2$ and so (2,4,1,2,1) holds, but $max~\hyperref[h18]{\phi_{20}}\leq10$, which gives contradiction.\vspace{1mm}\\
{\noindent \bf Case 127} $A>1, B>1, C>1, D\leq1, E\leq1, F\leq1, G\leq1, H\leq1, I\leq1, J>1$.\\
Proof. Claim(i)  $B>1.25$\\
For $B\leq1.25$, $max~\hyperref[h16]{\phi_{18}}\leq10$, which contradicts to (2,2,2,1,1,1,1).\\
Claim (ii)  $C<1.3$ \\
For $C\geq1.3$, $C^3/DEF>2$, so (2,4,2,1,1) holds, but $max~\hyperref[h13]{\phi_{15}}\leq10$\\
Claim(iii)  $A<1.62$ \\
For $A\geq1.62$, $A^3/BCD>2$ and $2G>H$, so (4,2,2,1,1) holds, but $max~\hyperref[h10]{\phi_{12}}\leq 10$.\\
Finally $B^3/CDE>2$ and so (1,4,2,1,1,1) holds. But $max~\hyperref[h14]{\phi_{16}}\leq10$.\vspace{1mm}\\
{\noindent \bf Case 189} $A>1, B>1, C\leq1, D>1, E\leq1, F\leq1, G\leq1, H\leq1, I>1, J>1$.\\
Proof. If $A^3/BCD\geq2$, then (4,2,1,1,1,1) holds, but $max~\hyperref[h15]{\phi_{17}}\leq10$, a contradiction.\\
If $A^3/BCD<2$, then $max~\hyperref[h20]{\phi_{22}}\leq10$, so $(1,2,3,1,1,1,1)$ contradicts.\vspace{1mm}\\
{\noindent \bf Case 249} $A>1, B>1, C\leq1, D\leq1, E\leq1, F\leq1, G\leq1, H>1, I>1, J>1$.\\
Proof.
Claim (i) $A^3/BCD<2$\\
If $A^3/BCD\geq2$, then (4,2,1,1,1,1) holds, but $max~\hyperref[h15]{\phi_{17}}\leq10$, a contradiction.\\
Claim (ii) $B^3/CDE<2$\\
If $B^3/CDE\geq2$, then (1,4,1,1,1,1) holds, but $max~\hyperref[h19]{\phi_{21}}\leq10$.\\
Finally $max~\hyperref[h16]{\phi_{18}}\leq10$, which contradicts to (2,2,2,1,1,1,1).\vspace{1mm}\\
{\noindent \bf Case 253} $A>1, B>1, C\leq1, D\leq1, E\leq1, F\leq1, G\leq1, H\leq1, I>1, J>1$.\\
Proof.
If $A^3/BCD\geq2$, then $(4,1,\cdots,1)$ holds, but $max~\hyperref[h22]{\phi_{24}}\leq10$.\\
If $A^3/BCD<2$, then $max~\hyperref[h16]{\phi_{18}}\leq10$, which contradicts to (2,2,2,1,1,1,1).\hfill$\Box$
\section{Most Difficult Cases}\label{Most Difficult}
In this section we  consider the remaining 7 cases, which are the most difficult to solve. These cases have also been solved  using the Optimization tool of the software Mathemetica, as in previous section. Here we use the obvious and common constraints, along with the 88 constraints which arise from all the possible \emph{strong} inequalities
corresponding to the partitions of 10 with summands equal to 1 and 2, provided the condition is satisfied for the particular part, e.g. we use $4G-2G^2/H$, only if $2G>H$. Few \emph{strong} inequalities are: \\
$4A-\frac{2A^2}{B}+C+D+E+F+G+H+I+J>10$,\\
$4A-\frac{2A^2}{B}+4C-\frac{2C^2}{D}+E+F+G+H+I+J>10$,\\
$A+4B-\frac{2B^2}{C}+4D-\frac{2D^2}{E}+4F-\frac{2F^2}{G}+H+I+J>10$,\\
$4A-\frac{2A^2}{B}+4C-\frac{2C^2}{D}+4E-\frac{2E^2}{F}+4G-\frac{2G^2}{H}+4I-\frac{2I^2}{J}>10$.
\begin{prop}
  The following cases do not arise:
\end{prop}
{\noindent \bf Case 16} $A > 1, B > 1, C > 1, D > 1, E > 1, F > 1, G \leq 1, H \leq 1, I\leq 1, J \leq 1$.\\
Proof. Here $2A>B, 2B>C, 2C>D, 2D>E, 2E>F, 2F>G, 2G>H, 2H>I, 2I>J$\\
Claim (i) $F<1.03$\\
If $F\geq1.03$, then $min~\hyperref[s6]{\psi_6}>2$, so the inequality (1,2,2,4,1) holds, but
$\max~\hyperref[h7]{\phi_9}\leq10.$ \\
Claim (ii) $E<1.202$\\
For $E\geq1.202$, $min~\hyperref[s5]{\psi_5}>2$, so the inequality (2,2,4,2) holds, but
$\max~\hyperref[h23]{\phi_{25}}\leq10$.\\
Claim (iii) $D<1.415$\\
For $D\geq1.415$, $min~\hyperref[s4]{\psi_4}>2$, so the inequality (2,1,4,2,1) holds, but
$max~\hyperref[h24]{\phi_{26}}\leq10$.\\
Claim(iv) $\frac{C^3}{DEF}<2$\\
Suppose $\frac{C^3}{DEF}\geq2$. Now if $\frac{G^3}{HIJ}\geq2$, then (2,4,4) holds, but $\max~\hyperref[h25]{\phi_{27}}<10$ and if $\frac{G^3}{HIJ}<2$, then $\max~\hyperref[h26]{\phi_{28}}\leq10$, which contradicts to (2,4,2,2).\\
Claim(v) $\frac{F^3}{GHI}<2$\\
If $\frac{F^3}{GHI}\geq2$, then the inequality (1,2,2,4,1) holds, but
$\max~\hyperref[h7]{\phi_9}<10$.\\
Finally
$\max~\hyperref[h9]{\phi_{11}}\leq10$, which contradicts to (1,2,2,2,2,1).\vspace{2mm}\\


{\noindent \bf Case 31} $A>1, B>1, C>1, D>1, E>1, F\leq1, G\leq1, H\leq1, I\leq1, J>1$.\\
Proof. Here $2A>B, 2B>C, 2C>D, 2D>E, 2E>F, 2F>G, 2G>H, 2H>I$\\
Claim(i) $E<1.03$  \\
For $E\geq1.03$, $\min~\hyperref[s5]{\psi_5}>2$, so the inequality (2,2,4,1,1) holds, but $\max~\hyperref[h5]{\phi_7}\leq10$.\\
Claim(ii)  $D<1.202$ \\
For $D\geq1.202$, $\min~\hyperref[s4]{\psi_4}>2$, so (1,2,4,2,1) holds, but $\max~\hyperref[h27]{\phi_{29}}\leq10$.\\
Claim(iii) $C<1.414$ \\
For $C\geq1.414$, $\min~\hyperref[s3]{\psi_3}>2$, so (2,4,2,1,1) holds, but $\max~\hyperref[h28]{\phi_{30}}\leq10$.\\
Claim(iv) $B<1.52$ \\
For $B\geq1.52$, $\min~\hyperref[s2]{\psi_2}>2$\\
If $\frac{F^3}{GHI}\geq2$, then (1,4,4,1) holds, but $\max~\hyperref[h8]{\phi_{10}}\leq10$.  \\
If $\frac{F^3}{GHI}<2$, then $\max~\hyperref[h8]{\phi_{13}}\leq10$, which contradicts to (1,4,2,2,1).\\
Claim (v) $A^3/BCD<2$\\ As for $A^3/BCD\geq2$, (4,2,2,1,1) holds, but $\max~\hyperref[h10]{\phi_{12}}\leq10$. \\
Claim (vi) $B>1.37$, $G<0.77$, $H<0.66$\\
As if $B\leq1.37$ or $G\geq0.77$ or $H\geq0.66$, the Max of $\max~\hyperref[h12]{\phi_{14}}\leq10$, which contradicts to (2,2,2,2,1,1)\\
Now working as in Claims(iii) and (iv) respectively we get $C<1.37$ and  $B^3/CDE<2$.\\
Finally we get $\frac{E^3}{FGH}>2$, so (2,2,4,1,1) holds , but $\max~\hyperref[h5]{\phi_{7}}\leq10.$\vspace{2mm}\\
{\noindent \bf Case 32} $A>1, B>1, C>1, D>1, E>1, F\leq1, G\leq1, H\leq1, I\leq1, J\leq1$.\\
Proof. Here $2A>B, 2B>C, 2C>D, 2D>E, 2E>F, 2F>G, 2G>H, 2H>I$\\
First take $2I\leq J$, then  $\frac{E^3}{FGH}>2$ and so the inequality (2,2,4,1,1) holds, but $\max~\hyperref[h5]{\phi_{7}}\leq10$.
So we must have $2I>J$, and consequently we can use $4I-2I^2/J$ instead of $2J$, wherever required.\\
 Now we have the following claims:\\
Claim (i) $E<1.208$\\
For $E\geq1.208$, $\min~\hyperref[s5]{\psi_{5}}>2$ and so the inequality (2,2,4,2) holds, but $\max~\hyperref[h23]{\phi_{25}}\leq10$.\\
Claim(ii) $D<1.418$\\
If $D\geq1.418$, then $\min~\hyperref[s4]{\psi_{4}}>2$ and so (2,1,4,2,1) holds, i.e. $4A-\frac{2A^2}{B}+C+4D-\frac{1}{2}\frac{D^4}{EFG}+4H-\frac{2H^2}{I}+J>10$. By taking this inequality as an additional constraint, we see that $\max~\hyperref[h27]{\phi_{29}}\leq10$, which contradicts to (1,2,4,2,1).\\
Claim (iii) $C<1.605$\\
For $C\geq 1.605$, $\min~\hyperref[s3]{\psi_{3}}>2$. Further if $G^3/HIJ\geq2$, the inequality (2,4,4) will hold, but $\max~\hyperref[h25]{\phi_{27}}\leq10$. If $G^3/HIJ<2$, then $\max~\hyperref[h26]{\phi_{28}}\leq10$, thus the inequality (2,4,2,2) contradicts.\\
Claim(iv) $B^3/CDE<2$\\
Suppose $B^3/CDE\geq2$. Further if $F^3/GHI\geq 2$, the inequality (1,4,4,1) holds, but $\max~\hyperref[h8]{\phi_{10}}\leq10$. If $F^3/GHI<2$, then $\max~\hyperref[h7]{\phi_{13}}\leq10$, thus the inequality (1,4,2,2,1) contradicts.\\
Claim(v) $A^3/BCD<2$ \\
Suppose $A^3/BCD\geq2$. Further if $E^3/FGH\geq 2$, the inequality (4,4,2) holds, but $\max~\hyperref[h29]{\phi_{31}}\leq10$. If $E^3/FGH<2$, then $\max~\hyperref[h30]{\phi_{32}}\leq10$, thus the inequality (4,2,2,2) contradicts.\\
Claim (vi) $B>1.38, H<0.775; I<0.7$\\
If either of  $B\leq 1.38$, $H\geq 0.775$, $I\geq0.7$ holds, then $\max~\hyperref[h9]{\phi_{11}}\leq10$, so the inequality (1,2,2,2,2,1) contradicts. \\
Claim (vii) $G<0.847$ and $C<1.555$\\
For $G\geq0.847$, then $\min~\hyperref[s7]{\psi_{7}}>2.$ Further if $C^3/DEF\geq2$, the inequality (2,4,4) will hold, but $\max~\hyperref[h25]{\phi_{27}}\leq10$. If $C^3/DEF<2$, then $\max~\hyperref[h1]{\phi_3}\leq10$, thus the inequality (2,2,2,4) contradicts.\\
Further working as in Claim(iii), we get $C<1.555$.\\
Claim(viii) $E<1.185$ and $D<1.394$\\
If $E\geq1.185$, then $E^3/FGH>2$, so the inequality (1,2,1,4,2) holds, i.e. $A+4B-2B^2/C+D+4E-(1/2)E^4/FGH+4I-2I^2/J>10$. By taking this inequality as an additional constraint, we see that $\max~\hyperref[h23]{\phi_{25}}\leq10$, which contradicts to (2,2,4,2).\\
Further working as in Claim(ii), we get $D<1.394$.\\
Claim(ix) $C\leq1.49$\\
Take $C>1.49$, then $\min~\hyperref[s3]{\psi_{3}}>2$. Further if $G^3/HIJ\geq2$, the inequality (2,4,4) contradict. So we must have $G^3/HIJ<2$, then working as in Claim(viii) we get $E<1.11$. Now $\max~\hyperref[h26]{\phi_{28}}\leq10$, thus the inequality (2,4,2,2) contradicts. So we must have $C\leq1.49$\\ Now working as in Claims (ii),(viii),(vi) respectively, we get $D<1.34$, $E<1.155$, $B>1.41$, $J<0.855$.
Claim(x) $F<0.989$\\
For $F\geq0.989$, $(F^3)/(GHI)>2$ and the inequality (1,2,2,4,1) holds, but $\max~\hyperref[h7]{\phi_{9}}\leq10$.\\
Further working as in Claims (vii),(ii),(viii) respectively we get $G<0.836$, $D<1.315$, $E<1.14$.\\
Claim(xii) $C<1.44$\\
Suppose $C^3/DEF>2$. Further if $G^3/HIJ\geq2$, the inequality (2,4,4) contradicts. So we must have $G^3/HIJ<2$, then working as in Claim(viii) we get $E<1.09$, which contradicts to (2,4,2,2).\\
Finally working as in Claim(ii), we get $D^3/EFG<2$, which further gives $\min~\hyperref[s5]{\psi_{5}}>2$, but $\max~\hyperref[h23]{\phi_{25}}\leq10$, which contradicts to (2,2,4,2).\vspace{2mm}\\
{\noindent \bf Case 61} $A>1, B>1, C>1, D>1, E\leq1, F\leq1, G\leq1, H\leq1, I>1, J>1$.\\
Proof. Here $2A>B, 2B>C, 2C>D, 2D>E, 2E>F, 2F>G, 2G>H$\\
Claim(i) $D<1.03$\\
For $D\geq1.03$, $\min~\hyperref[s4]{\psi_4}>2$, and so the inequality (1,2,4,1,1,1) holds, but $\max~\hyperref[h31]{\phi_{33}}\leq10$.\\
Claim(ii)  $C<1.2$\\
For $C\geq1.2$, $\min~\hyperref[s3]{\psi_3}>2$ and so the inequality (2,4,2,1,1) holds, but $\max~\hyperref[h28]{\phi_{30}}\leq10$.\\
Claim(iii) $B<1.414$\\
For $B\geq1.414$, $\min~\hyperref[s2]{\psi_2}>2$ and so (1,4,2,1,1,1) holds, but $\max~\hyperref[h14]{\phi_{16}}\leq10$.\\
Claim(iv) $A^3/BCD<2$\\
As for $A^3/BCD\geq2$, the inequality (4,2,2,1,1) holds, but $\max~\hyperref[h10]{\phi_{12}}\leq10$. \\
Claim(v) $G<0.66$\\
For $G\geq0.66$, $\max~\hyperref[h12]{\phi_{14}}\leq10$, which contradicts to the inequality (2,2,2,2,1,1).\\
Finally $\min~\hyperref[s4]{\psi_4}>2$ and we get contradiction working as in Claim(i).\vspace{2mm}\\

{\noindent \bf Case 63} $A>1, B>1, C>1, D>1, E\leq1, F\leq1, G\leq1, H\leq1, I\leq1, J>1$.\\
Proof. Here $2A>B, 2B>C, 2C>D, 2D>E, 2E>F, 2F>G, 2G>H$\\
IF $H<0.5$ then $\min~\hyperref[s3]{\psi_3}>2$, so the inequality (2,4,2,1,1) holds, but $\max~\hyperref[h13]{\phi_{15}}\leq10$, which contradicts .\\
Thus we have $H\geq0.5$, then $2H>I$, and consequently we can use $4H-2(H^2)/I$ instead of $2I$, wherever required. Now we have following claims:\\
Claim(i) $D<1.2$\\
If $D\geq1.2$, then $\min~\hyperref[s4]{\psi_4}>$ and so the inequalities  (2,1,4,2,1), (3,4,2,1) hold, i.e. $4A-2A^2/B+C+4D-(1/2)D^4/EFG+4H-2H^2/I+J>10$ and $4A-A^3/BC+4D-(1/2)D^4/EFG+4H-2H^2/I+J>10$. But by taking these two inequalities as additional constraints, we find that $\max~\hyperref[h27]{\phi_{29}}\leq10$, which contradicts to (1,2,4,2,1).\\
Claim(ii) $C<1.418$\\
For $C\geq1.418,$ $\min~\hyperref[s3]{\psi_3}>2$, so the inequality (2,4,2,1,1) holds, but $\max~\hyperref[h13]{\phi_{15}}\leq10$.\\
Claim(iii)  $B<1.605$\\
For $B\geq1.605$, $\min~\hyperref[s2]{\psi_2}>2$. If $F^3/GHI>2$, then the inequality (1,4,4,1) holds, but $\max~\hyperref[h8]{\phi_{10}}\leq10.$
If $F^3/GHI\leq2$, the $\max~\hyperref[h11]{\phi_{13}}\leq10$, which contradicts to (1,4,2,2,1).\\
Claim(iv) $A^3/BCD<2$\\
Suppose $A^3/BCD\geq2$.\\
If $E^3/FGH\geq2$, then (4,4,1,1) holds, but $\max~\hyperref[h32]{\phi_{34}}\leq10$.\\
If $E^3/FGH<2$, the $\max~\hyperref[h10]{\phi_{12}}\leq10$, which contradicts to (4,2,2,1,1).\\
Claim(v) $G<0.78$ and $B<1.555$  \\
As for $G\geq0.78$, we get $\max~\hyperref[h12]{\phi_{14}}\leq10$, contradicting (2,2,2,2,1,1).\\
Further working as in Claim (iii), we get $B<1.555$.\\
Claim(vi) $B<1.51$\\
Take $B>1.51$, then $\min~\hyperref[s2]{\psi_2}>2$. Further if $F^3/GHI\geq2$, we get contradiction to (1,4,4,1). So $F^3/GHI<2$. Now working as in Claim(i), we get $D<1.11$. This further implies that   $\max~\hyperref[h11]{\phi_{13}}\leq10$, which contradicts to (1,4,2,2,1).\\
Now working as in Claims(ii) and (i) respectively, we get $C<1.32$ and $D<1.145$.\\
Claim(vii) $F<0.837$\\
For $F\geq0.837$, $F^3/GHI>2$, so the inequality (1,2,2,4,1) holds, but $\max~\hyperref[h7]{\phi_{9}}\leq10$.\\
Claim(viii) $A>1.42$ and  $B>1.39$\\
If $A\leq1.42$ or  $B\leq1.39$, we get $\max~\hyperref[h12]{\phi_{14}}\leq10$, contradicting (2,2,2,2,1,1).\\
Now working as in Claims(ii), (i) and (iii) respectively, we get $C<1.307$, $D<1.133$, $B<1.49$\\
Claim(ix) $E<0.975$\\
For $E\geq0.975$, $E^3/FGH>2$, so the inequality (2,2,4,1,1) holds, but $\max~\hyperref[h33]{\phi_{35}}\leq10$.\\
Now working as in Claim (ii), we get contradiction if $C^3/DEF\geq2$.\\
So $C^3/DEF<2$. Finally we get $\min~\hyperref[s2]{\psi_2}>2$ and $\max~\hyperref[h11]{\phi_{13}}\leq10$, which contradicts to (1,4,2,2,1).\vspace{2mm}\\
{\noindent \bf Case 64} $A>1, B>1, C>1, D>1, E\leq1, F\leq1, G\leq1, H\leq1, I\leq1, J\leq1$.\\
Proof. Here $2A>B, 2B>C, 2C>D, 2E>F, 2F>G, 2G>H$\\
Claim(i) $A^3/BCD<2$\\
Suppose $A^3/BCD\geq2$, then the inequality (4,2,2,1,1) holds, but $\max~\hyperref[h10]{\phi_{12}}\leq10$.\\
Claim(ii) $D<1.268$\\
If $D\geq1.268$, then (3,4,2,1) holds, but $\max~\hyperref[h34]{\phi_{36}}\leq10$.\\
Claim(iii) $C^3/DEF<2$\\
For $C^3/DEF\geq2$, (2,4,2,1,1) holds, but $\max~\hyperref[h28]{\phi_{30}}\leq10$.\\
Finally $\max~\hyperref[h12]{\phi_{14}}\leq10$, which contradicts to (2,2,2,2,1,1).\vspace{2mm}\\
{\noindent \bf Case 125} $A>1, B>1, C>1, D\leq1, E\leq1, F\leq1, G\leq1, H\leq1, I>1, J>1$.\\
Proof. Here $2A>B, 2B>C, 2C>D, 2D>E, 2E>F, 2F>G$\\
First suppose $2G<H$, then $\max~\hyperref[h35]{\phi_{37}}\leq10$, which contradicts to (1,2,2,2,1,1,1).\\
Now suppose $2G>H$. We have the following claims:\\
Claim(i) $C<1.2$\\
If $C\geq1.2$, then $\min~\hyperref[s3]{\psi_3}>2$. So (2,4,2,1,1) holds, but $\max~\hyperref[h13]{\phi_{15}}\leq10$.\\
Claim(ii) $B<1.421$.
For $B\geq1.421$, $\min~\hyperref[s2]{\psi_2}>2$, so (1,4,2,1,1,1) holds, but $\max~\hyperref[h14]{\phi_{16}}\leq10$.\\
Claim(iii) $A^3/BCD<2$\\
For $A^3/BCD\geq2$, (4,2,1,1,1,1) holds, but $\max~\hyperref[h15]{\phi_{17}}\leq10$.\\
Claim(iv) $F<0.79, G<0.7, H<0.84$\\
If $F\geq0.79 ~or~ G\geq0.7 ~or~ H\geq0.84$, then $\max~\hyperref[h35]{\phi_{37}}\leq10$, which contradicts to (1,2,2,2,1,1,1).\\
 Claim(v) $E<0.85$\\
 For $E\geq0.85,$ $\min~\hyperref[s5]{\psi_5}>2$, so (2,2,4,1,1) holds. But $\max~\hyperref[h5]{\phi_{7}}\leq10.$\\
Finally $\min~\hyperref[s3]{\psi_3}>2$ and we get contradiction as in claim(i) to the inequality (2,4,2,1,1).\hfill$\Box$\vspace{2mm}\\
This settles all the cases and hence completes the proof of Woods' Conjecture for $n=10$.

\end{document}